**Title: A Game Theoretic Approach to Sustainizability® Over Sets and its application to a multi-species population model**

**Order of Authors:** Ioannis V. Manousiouthakis[a], Vasilios I. Manousiouthakis[b*]
**First Author:** Ioannis V. Manousiouthakis[a], [a] Hydrogen Engineering Research Company, LLC (H-E-R-C, LLC), Los Angeles, CA, 90077, United States of America
**Corresponding Author:** Vasilios I. Manousiouthakis[b], [b] Department of Chemical and Biomolecular Engineering, Hydrogen Engineering Research Consortium (HERC), University of California at Los Angeles (UCLA), Los Angeles, CA, USA (ORCID: 0000-0002-5926-9923)
Email: vasilios@ucla.edu, Phone: (310) 206-0300, Address: 5549 Boelter Hall, Box 951592, Los Angeles, CA, 90095-1592, U.S.A.

**Abstract:** In this work, a theorem is first proved which presents a game theoretic formulation of a necessary and sufficient sustainizability® over a set condition for a general system described by ordinary differential equations (ODEs). Then, two additional theorems are proved for the n-species Gause-Lotka-Volterra (GLV) population model, establishing necessary and sufficient sustainability and sustainizability® conditions over rectangular sets. Three case studies on the May-Leonard, 3-species, GLV model are then presented to illustrate the power of the above Theorems. In two of these case studies (2 and 3), it is shown that a particular instance of the 3-species, GLV model is unsustainable but sustainizable through allowable action.

## Introduction

The recent focus of the United Nations on sustainable development, crystallized in the statement "Humanity has the ability to make development sustainable to ensure that it meets the needs of the present without compromising the ability of future generations to meet their own needs." contained in the report "Our Common Future" (Brundtland, U.N. World Commission, 1987), has driven in recent years a dramatic increase in research efforts to define and assess sustainability. In 1992, Rees identified (Rees, 1992) "The total area of land required to sustain an urban region" as its "ecological footprint". The incorporation of ecological considerations in the design of manufacturing systems has also been considered (Bakshi and Fiksel, 2003), while the Sustainability Assessment by Fuzzy Evaluation (SAFE) hierarchical fuzzy inference conceptual framework (Phillis and coworkers, 2001, 2003, 2004, 2009) crystallizes both ecological and human considerations in the form of a sustainability index. The Sustainability Interval Index concept (Conner et al, 2012) was subsequently introduced as another sustainability assessment tool that employs both fuzzy logic and interval analysis to account for uncertainty in sustainable performance basic indicator data. The concept of categorizing sustainability assessment models as a system of systems (Phillis et al, 2010), has also been put forward as a means of accounting for system complexity and for ecological and human aspects of sustainable system performance. Over the last few decades the U.S. EPA has also been focusing on sustainability. In addressing the question "What is Sustainability Anyway?", the U.S. EPA's National Risk Management Research Laboratory came up with the "straw man" definition ''sustainability occurs when we maintain or improve the material and social conditions for human health and the environment over time without exceeding the ecological capabilities that support them.'' (Sikdar, 2003). An additional sustainability hypothesis originating from the same organization stated ''sustainable systems do not lose or gain Fisher information over time.'' (Cabezas and Fath, 2003). The EPA has also focused on developing the sustainability assessment tool GREENSCOPE (Gauging Reaction Effectiveness for the ENvironmental Sustainability of Chemistries with a Multi-Objective Process Evaluator) to assess and improve the design of chemical processes, in regard to several Material Efficiency, Energy, Economic, and Environmental, sustainable performance indicators (Gonzalez and Smith, 2003).

Two recently proposed sustainability concepts that are the focus of this work, are the Sustainability Over Sets (SOS) and Sustainizability® Over Sets (SIZOS) (Manousiouthakis and coworkers, 2018, 2019, 2020, 2021), which are mathematically rigorous and computationally efficient methods for sustainable system analysis and synthesis respectively. Both methods first define a set in the considered system's state space that can reflect input from discipline experts as to whether the system is deemed sustainable if its dynamic evolution trajectory remains within the set. SOS identifies whether the system's trajectories initiated within the considered set remain for all time within the set, while SIZOS identifies whether there exist admissible control strategies that can ensure the system's trajectories initiated within the considered set remain for all time within the set. Both methods, employ the concept of positive invariant sets in carrying out sustainable system analysis and synthesis. The computation of Positively Invariant Sets has been pursued through various methods, such as Newton's method (Baier et al, 2010), and variational techniques (Junge and Kevrekidis, 2017).

In this work, the SOS and SIZOS methods are first briefly reviewed, then a novel game theoretic approach to the sustainizability® over a set for a general system is presented in a Theorem, which yields a necessary and sufficient sustainizability® condition, and then SOS and SIZOS are applied to the n-species Gause-Lotka-Volterra (GLV) population model, for which two Theorems are proved and three case studies on the 3-species, GLV model are carried out.

**Sustainability Over Sets (SOS) Analysis**

Sustainability Over Sets (SOS) analysis for a time-invariant unforced system, modeled as an initial value problem for a system of ordinary differential equations (ODEs), is carried out by identifying whether the system's trajectories initiated within a considered set remain for all time within the set, i.e. by identifying whether the considered set is a Positively Invariant Set for the system. The computation of Positively Invariant Sets has been pursued through various methods, such as Newton's method (Baier et al, 2010), and variational techniques (Junge and Kevrekidis, 2017). Mathematically, SOS is stated as follows:

Let the system

$$\frac{dx_i(t)}{dt} = f_i\left(\left\{x_j(t)\right\}_{j=1}^{n}\right)\ \forall t \in [0,\infty)\ \ x_i(0) = x_i^0\ \ i \in \{1,\ldots,n\} \tag{1}$$

admit a unique solution $\left\{x_i^1\left(t,\left\{x_j^0\right\}_{j=1}^{n}\right)\right\}_{i=1}^{n}\ \ \forall t \in [0,\infty)\ \forall \left\{x_j^0\right\}_{j=1}^{n} \in \mathbb{R}^n$ with no finite time escape.

Let $\Phi_k : \mathbb{R}^n \to \mathbb{R}\ \ \forall k \in \{1,\ldots,m\}$ be continuously differentiable functions on $\mathbb{R}^n$, and define the uncountably infinite collection of index sets

$$S\left(\left\{z_j\right\}_{j=1}^{n}\right) \triangleq \left\{k \in \{1,\ldots,m\} : \Phi_k\left(\left\{z_j\right\}_{j=1}^{n}\right) = 0\right\}\ \ \forall \left\{z_j\right\}_{j=1}^{n} \in \mathbb{R}^n \tag{2}$$

Consider the closed set $F$, with complement $F^c$, and boundary $\partial F$ defined as:

$$\begin{gathered} F \triangleq \left\{\left\{z_j\right\}_{j=1}^{n} \in \mathbb{R}^n : \Phi_k\left(\left\{z_j\right\}_{j=1}^{n}\right) \le 0\, \forall k \in \{1,\ldots,m\}\right\}, F^c \triangleq \left\{\left\{z_j\right\}_{j=1}^{n} \in \mathbb{R}^n : \exists k \in \{1,\ldots,m\} : \Phi_k\left(\left\{z_j\right\}_{j=1}^{n}\right) > 0\right\} \\ \partial F \triangleq F \cap \overline{F^c} \triangleq \left\{\left\{z_j\right\}_{j=1}^{n} \in \mathbb{R}^n : \Phi_k\left(\left\{z_j\right\}_{j=1}^{n}\right) \le 0\ \forall k \in \{1,\ldots,m\} \wedge S\left(\left\{z_j\right\}_{j=1}^{n}\right) \ne \varnothing\right\} \end{gathered} \tag{3}$$

<u>The unforced system (1) is sustainable over the set $F \subset \mathbb{R}^n$ (SOS $F \subset \mathbb{R}^n$)</u> if and only if:

$$\left\{\left\{x_i^1\left(t,\left\{x_j^0\right\}_{j=1}^n\right)\right\}_{i=1}^n \in F \quad \forall t\in[0,\infty) \ \ \forall\left\{x_j^0\right\}_{j=1}^n \in F\right\} \Leftrightarrow$$

$$\left\{\forall\left\{z_j\right\}_{j=1}^n \in F:\ S\left(\left\{z_j\right\}_{j=1}^n\right)\neq\varnothing, \left[\sum_{i=1}^n \frac{\partial\Phi_k\left(\left\{z_j\right\}_{j=1}^n\right)}{\partial z_i}\cdot f_i\left(\left\{z_j\right\}_{j=1}^n\right)\right]\leq 0 \ \forall k\in S\left(\left\{z_j\right\}_{j=1}^n\right)\right\} \tag{4}$$

For the special case of a closed rectangular set, i.e.:

$$F^R \triangleq \left\{\left\{z_j\right\}_{j=1}^n \in \mathbb{R}^n : \left\{z_j^l \leq z_j \leq z_j^u,\ z_j^l < z_j^u\right\} \ \ \forall j\in\{1,\ldots,n\}\right\} \tag{5}$$

the considered ODE system is <u>SOS</u> $F^R$ iff:

$$\left\{\begin{array}{l} \left\{f_i\left(\left\{x_j^{i,u}\right\}_{j=1}^n\right)\leq 0 \ \forall\left\{x_j^{i,u}\right\}_{j=1}^n : \left[\begin{array}{l} x_i^{i,u}=z_i^u \\ z_j^l \leq x_j^{i,u} \leq z_j^u \ \ \forall j\in\{1,\ldots,n\};\ j\neq i\end{array}\right]\right\} \ \forall i\in\{1,\ldots,n\} \\ \left\{f_i\left(\left\{x_j^{i,l}\right\}_{j=1}^n\right)\geq 0 \ \forall\left\{x_j^{i,l}\right\}_{j=1}^n : \left[\begin{array}{l} x_i^{i,l}=z_i^l \\ z_j^l \leq x_j^{i,l} \leq z_j^u \ \ \forall j\in\{1,\ldots,n\};\ j\neq i\end{array}\right]\right\} \ \forall i\in\{1,\ldots,n\}\end{array}\right\} \tag{6}$$

**Sustainizability® Over Sets (SIZOS) Synthesis**

Sustainizability® Over Sets (SIZOS) synthesis for a time-invariant forced system, modeled as an initial value problem for a system of ordinary differential equations (ODEs), is carried out by identifying whether there exist admissible control strategies that can ensure the system's trajectories initiated within a considered set remain for all time within the set (Jorat and Manousiouthakis, 2019). Mathematically, this can be stated as:

Let the system

$$\frac{dx_i(t)}{dt}=f_i\left(\left\{x_j(t)\right\}_{j=1}^n,\left\{u_l(t)\right\}_{l=1}^p\right),\ \left\{u_l(t)\right\}_{l=1}^p \in U\subset\mathbb{R}^p \ \ \forall t\in[0,\infty), x_i(0)=x_i^0 \ \ i\in\{1,\ldots,n\} \tag{7}$$

admit the unique solution

$$\left\{x_i^2\left(t,\left\{u_l(t)\right\}_{l=1}^p,\left\{x_j^0\right\}_{j=1}^n\right)\right\}_{i=1}^n \quad \forall t\in[0,\infty) \ \ \forall\left\{x_j^0\right\}_{j=1}^n\in\mathbb{R}^n \ \ \forall\left\{u_l\right\}_{l=1}^p : \left\{\begin{array}{l}[0,\infty)\rightarrow U\subset\mathbb{R}^p \\ t\rightarrow\left\{u_l(t)\right\}_{l=1}^p \in U\subset\mathbb{R}^p\end{array}\right\}$$

with no finite time escape behavior. Using state feedback strategies

$$\left\{\left\{\bar{u}_l\right\}_{l=1}^p : \mathbb{R}^n \rightarrow U\subset\mathbb{R}^p,\ \left\{\bar{u}_l\right\}_{l=1}^p : \left\{x_j\right\}_{j=1}^n \rightarrow \left\{\bar{u}_l\left(\left\{x_j\right\}_{j=1}^n\right)\right\}_{l=1}^p \in U\subset\mathbb{R}^p\right\} \tag{8}$$

the above forced system is brought into the following unforced system form:

$$\left\{\begin{array}{l} \dfrac{dx_i(t)}{dt}=f_i\left(\left\{x_j(t)\right\}_{j=1}^n,\left\{\bar{u}_l\left(\left\{x_j(t)\right\}_{j=1}^n\right)\right\}_{l=1}^p\right)\triangleq g_i\left(\left\{x_j(t)\right\}_{j=1}^n\right) \\ \left\{\bar{u}_l\left(\left\{x_j(t)\right\}_{j=1}^n\right)\right\}_{l=1}^p \in U\subset\mathbb{R}^p,\ x_i(0)=x_i^0 \ \ i\in\{1,\ldots,n\}\end{array}\right\} \ \forall t\in[0,\infty) \tag{9}$$

whose unique solution $\left\{x_i^3\left(t,\left\{x_j^0\right\}_{j=1}^n,\left\{\bar{u}_l\left(\left\{x_j^3(t)\right\}_{j=1}^n\right)\right\}_{l=1}^p\right)\right\}_{i=1}^n \quad \forall t\in[0,\infty) \ \ \forall\left\{x_j^0\right\}_{j=1}^n\in\mathbb{R}^n$

exhibits no finite time escape behavior.

Then the SIZOS $F$, $U$ definition and associated necessary and sufficient conditions are:

$$\left\{\begin{array}{l}\left\{\begin{array}{l}\dfrac{dx_i(t)}{dt}=f_i\left(\left\{x_j(t)\right\}_{j=1}^{n},\left\{u_l(t)\right\}_{l=1}^{p}\right)\ \forall t\in[0,\infty)\ i\in\{1,\ldots,n\}\\ \left\{u_l(t)\right\}_{l=1}^{p}\in U\subset\mathbb{R}^p\ \ \forall t\in[0,\infty),\ x_i(0)=x_i^0\ i\in\{1,\ldots,n\}\end{array}\right\}\\ \textit{is sustainizable over the sets } F\subset\mathbb{R}^n,\ U\subset\mathbb{R}^p\end{array}\right\}\Leftrightarrow$$

$$\left\{\exists\left\{u_l\right\}_{l=1}^{p}:\left\{\begin{array}{l}[0,\infty)\to U\\ t\to\left\{u_l(t)\right\}_{l=1}^{p}\in U\end{array}\right\}:\left\{x_i^2\left(t,\left\{u_l(t)\right\}_{l=1}^{p},\left\{x_j^0\right\}_{j=1}^{n}\right)\right\}_{i=1}^{n}\in F\left\{\begin{array}{l}\forall t\in[0,\infty)\\ \forall\left\{x_j^0\right\}_{j=1}^{n}\in F\end{array}\right\}\right\}\Leftrightarrow$$

$$\left\{\begin{array}{l}\left\{\exists\left\{\bar{u}_l\right\}_{l=1}^{p}:\mathbb{R}^n\to U\subset\mathbb{R}^p,\exists\left\{\bar{u}_l\right\}_{l=1}^{p}:\left\{x_j(t)\right\}_{j=1}^{n}\to\left\{\bar{u}_l\left(\left\{x_j(t)\right\}_{j=1}^{n}\right)\right\}_{l=1}^{p}\in U\right\}:\\ \left\{\begin{array}{l}\dfrac{dx_i(t)}{dt}=f_i\left(\left\{x_j(t)\right\}_{j=1}^{n},\left\{\bar{u}_l\left(\left\{x_j(t)\right\}_{j=1}^{n}\right)\right\}_{l=1}^{p}\right)\triangleq g_i\left(\left\{x_j(t)\right\}_{j=1}^{n}\right)\ \left\{\begin{array}{l}\forall t\in[0,\infty)\\ x_i(0)=x_i^0\\ \forall i\in\{1,\ldots,n\}\end{array}\right\}\\ \textit{is sustainable over the set } F\subset\mathbb{R}^n\end{array}\right\}\end{array}\right\}\Leftrightarrow$$

$$\left\{\exists\left\{\bar{u}_l\right\}_{l=1}^{p}:\left\{\begin{array}{l}F\to U\\ \left\{x_j\right\}_{j=1}^{n}\to\left\{\bar{u}_l\left(\left\{x_j\right\}_{j=1}^{n}\right)\right\}_{l=1}^{p}\in U\end{array}\right\}:\left\{x_i^3\left(\begin{array}{l}t,\left\{x_j^0\right\}_{j=1}^{n}\\ \left\{\bar{u}_l\left(\left\{x_j^3(t)\right\}_{j=1}^{n}\right)\right\}_{l=1}^{p}\end{array}\right)\right\}_{i=1}^{n}\in F\left\{\begin{array}{l}\forall t\in[0,\infty)\\ \forall\left\{x_j^0\right\}_{j=1}^{n}\in F\end{array}\right\}\right\}\Leftrightarrow$$

$$\left\{\forall\left\{z_j\right\}_{j=1}^{n}:\left\{\begin{array}{l}\left\{z_j\right\}_{j=1}^{n}\in F\\ S\left(\left\{z_j\right\}_{j=1}^{n}\right)\neq\varnothing\end{array}\right\}\exists\left\{\bar{u}_l\right\}_{l=1}^{p}:\left\{\begin{array}{l}\left\{\bar{u}_l\right\}_{l=1}^{p}:F\to U,\ \left\{\bar{u}_l\right\}_{l=1}^{p}:\left\{x_j\right\}_{j=1}^{n}\to\left\{\bar{u}_l\left(\left\{x_j\right\}_{j=1}^{n}\right)\right\}_{l=1}^{p}\in U\\ \sum_{i=1}^{n}\left[\begin{array}{l}\dfrac{\partial\Phi_k\left(\left\{z_j\right\}_{j=1}^{n}\right)}{\partial z_i}\cdot\\ f_i\left(\left\{z_j\right\}_{j=1}^{n},\left\{\bar{u}_l\left(\left\{z_j\right\}_{j=1}^{n}\right)\right\}_{l=1}^{p}\right)\end{array}\right]\leq 0\,\forall k\in S\left(\left\{z_j\right\}_{j=1}^{n}\right)\end{array}\right\}\right\}\quad(10)$$

Then, a game theoretic necessary and sufficient condition for Sustainizability® over a set is:

<u>Theorem 1</u>

Consider the system

$\dfrac{dx_i(t)}{dt}=f_i\left(\left\{x_j(t)\right\}_{j=1}^{n},\left\{u_l(t)\right\}_{l=1}^{p}\right),\ \left\{u_l(t)\right\}_{l=1}^{p}\in U\subset\mathbb{R}^p\ \forall t\in[0,\infty),x_i(0)=x_i^0\ i\in\{1,\ldots,n\}$ that has no finite time escape behavior and admits the unique solution

$$\left\{x_i^2\left(t,\left\{u_l(t)\right\}_{l=1}^{p},\left\{x_j^0\right\}_{j=1}^{n}\right)\right\}_{i=1}^{n}\ \ \forall t\in[0,\infty)\ \ \forall\left\{x_j^0\right\}_{j=1}^{n}\in\mathbb{R}^n\ \ \forall\left\{u_l\right\}_{l=1}^{p}:\left\{\begin{array}{l}[0,\infty)\to U\subset\mathbb{R}^p\\ t\to\left\{u_l(t)\right\}_{l=1}^{p}\in U\subset\mathbb{R}^p\end{array}\right\}$$

Then the system is sustainizable over the sets $F\subset\mathbb{R}^n,\ U\subset\mathbb{R}^p$ (SIZOS $F,U$) if and only if

$$\left\{ 0 \geq \max_{\substack{\{z_j\}_{j=1}^n \in F \\ S\left(\{z_j\}_{j=1}^n\right) \neq \varnothing}} \min_{\{\bar{\bar{u}}_l\}_{l=1}^p \in U} \max_{k \in S\left(\{z_j\}_{j=1}^n\right)} \sum_{i=1}^n \left[ \frac{\partial \Phi_k\left(\{z_j\}_{j=1}^n\right)}{\partial z_i} \cdot f_i\left(\{z_j\}_{j=1}^n, \{\bar{\bar{u}}_l\}_{l=1}^p\right) \right] \forall \{z_j\}_{j=1}^n : \left\{ \begin{matrix} \{z_j\}_{j=1}^n \in F \\ S\left(\{z_j\}_{j=1}^n\right) \neq \varnothing \end{matrix} \right\} \right\} \tag{11}$$

For the special case of the closed rectangular set $F^R \triangleq \left\{ \{z_j\}_{j=1}^n \in \mathbb{R}^n : \{z_j^l \leq z_j \leq z_j^u,\, z_j^l < z_j^u\} \;\; \forall j \in \{1,\ldots,n\} \right\}$ defined in (5), define the sets

$$\left\{ \begin{matrix} S^U\left(\{z_j\}_{j=1}^n\right) \triangleq \{k \in \{1,\ldots,n\} : z_k = z_k^u\} \\ S^L\left(\{z_j\}_{j=1}^n\right) \triangleq \{k \in \{1,\ldots,n\} : z_k = z_k^l\} \end{matrix} \right\}, \;\; S^U\left(\{z_j\}_{j=1}^n\right) \cap S^L\left(\{z_j\}_{j=1}^n\right) = \varnothing \tag{12}$$

Then, the considered forced system (7) is (SIZOS $F^R \subset \mathbb{R}^n, U \subset \mathbb{R}^p$) iff:

$$\left\{ \left\{ \forall \{x_j\}_{j=1}^n : \left[ \begin{matrix} \left\{ \begin{matrix} S^U\left(\{x_j\}_{j=1}^n\right) \neq \varnothing \\ \vee \\ S^L\left(\{x_j\}_{j=1}^n\right) \neq \varnothing \end{matrix} \right\} \\ x_j^l \leq x_j \leq x_j^u \;\; \forall j \in \{1,\ldots,n\} \end{matrix} \right] \right\} \exists \{\bar{u}_l\}_{l=1}^p : \left\{ \begin{matrix} \{\bar{u}_l\}_{l=1}^p : \{x_j\}_{j=1}^n \to \left\{ \bar{u}_l\left(\{x_j\}_{j=1}^n\right) \right\}_{l=1}^p \in U \subset R^p \\ f_i\left( \{x_j\}_{j=1}^n, \left\{ \bar{u}_l\left(\{x_j\}_{j=1}^n\right) \right\}_{l=1}^p \right) \leq 0 \;\; \forall i \in S^U\left(\{x_j\}_{j=1}^n\right) \\ -f_i\left( \{x_j\}_{j=1}^n, \left\{ \bar{u}_l\left(\{x_j\}_{j=1}^n\right) \right\}_{l=1}^p \right) \leq 0 \;\; \forall i \in S^L\left(\{x_j\}_{j=1}^n\right) \end{matrix} \right\} \right\} \Leftrightarrow$$

$$\left\{ 0 \geq \max_{\substack{x_j^l \leq x_j \leq x_j^u \;\; \forall j \in \{1,\ldots,n\} \\ \left\{ \begin{matrix} S^U\left(\{x_j\}_{j=1}^n\right) \neq \varnothing \\ \vee \\ S^L\left(\{x_j\}_{j=1}^n\right) \neq \varnothing \end{matrix} \right\}}} \min_{\{\bar{\bar{u}}_l\}_{l=1}^p \in U} \max \left[ \begin{matrix} \max_{k \in S^U\left(\{x_j\}_{j=1}^n\right)} \left[ f_i\left( \{x_j\}_{j=1}^n, \{\bar{\bar{u}}_l\}_{l=1}^p \right) \right], \\ \max_{k \in S^L\left(\{x_j\}_{j=1}^n\right)} \left[ -f_i\left( \{x_j\}_{j=1}^n, \{\bar{\bar{u}}_l\}_{l=1}^p \right) \right] \end{matrix} \right] \right\} \tag{13}$$

Having outlined and quantified the Sustainability Over Sets (SOS) and Sustainizability® Over Sets (SIZOS) concepts, SOS analysis and SIZOS synthesis over a rectangular set is next carried out for the considered multispecies GLV population model, and two Theorems are proved.

**Sustainability and Sustainizability® of the n-species Gause-Lotka-Volterra GLV model**

Lotka developed in 1910 (Lotka, 1910) a model for the behavior of autocatalytic chemical reactions, that he subsequently extended in 1925 (Lotka, 1925) to models capturing the behavior of a system consisting of a host and a parasite, and a biological system consisting of three species, with the first and third species feeding on a constant resource, and being fed by the second species. In 1926, Volterra considered (Volterra, 1926) a case of two species "of which one, finding sufficient food in its environment, would multiply indefinitely when left to itself, while the other would perish for lack of nourishment if left alone; but the second feeds upon the first, and so the two species can co-exist together." He subsequently stated: "it is possible to establish two differential equations of the first order, non-linear, which can be integrated." He then went on to discuss a three species case. In 1934, Gause carried out (Gause, 1934) "an experimental investigation of the processes of the struggle for existence among unicellular organisms. Experiments on the competition between two species for a common place in the microcosm agreed completely with Volterra's theoretical equations, but as regards the processes of one species devouring another our results are not concordant with the forecasts of the mathematical theory." Later in the same manuscript he provided an explanation for theory-experiment discrepancies near population elimination conditions, stating: "This showed that when the number of individuals becomes reduced, and the conditions in the microcosm complicated, instead of the "deterministic" processes subject to differential equations we are confronted with "probabilities of change" in one direction or another."

Following these early developments, the Gause-Lotka-Volterra (GLV) population models have been the subject of many studies focusing on understanding the dynamics of natural populations of predators and prey.

In the original study of the nonlinear aspects of competition among three species (May and Leonard, 1975), the following general n-species GLV population model is considered:

$$\frac{dN_i(t)}{dt} = r_i N_i(t)\left[1 - \sum_{j=1}^{n} \alpha_{ij} N_j(t)\right] \ \forall i \in \{1,\ldots,n\}. \tag{14}$$

The above system has a multitude of equilibrium points, which considering that $r_i \neq 0 \ \forall i = 1,n$ must satisfy the equations

$$0 = N_{i,e}\left[1 - \sum_{j=1}^{n} \alpha_{ij} N_{j,e}\right] \ \forall i \in \{1,\ldots,n\}. \tag{15}$$

This set of nonlinear equations has $2^n$ equilibrium points, some of which may not be physically realizable, and only one of which may not involve the extinction of at least one of the species. To streamline the analysis that follows, the following index sets are first defined:

$$\left\{\begin{array}{l} \left\{ \mathrm{A}_i^+ \triangleq \left\{ j \in \{1,\ldots,n\} \ : \alpha_{ij} > 0 \right\}, \mathrm{A}_i^- \triangleq \left\{ j \in \{1,\ldots,n\} \ : \alpha_{ij} < 0 \right\}, \mathrm{A}_i^+ \cap \mathrm{A}_i^- = \varnothing \right\} \ \forall i \in \{1,\ldots,n\} \\ R^+ \triangleq \left\{ i \in \{1,\ldots,n\} \ : r_i > 0 \right\},\ R^- \triangleq \left\{ i \in \{1,\ldots,n\} \ : r_i < 0 \right\},\ R^+ \cup R^- = \{1,\ldots,n\},\ R^+ \cap R^- = \varnothing \end{array}\right\} \tag{16}$$

<u>Theorem 2</u>

Consider the rectangular set

$$F^{RN} \triangleq \left\{ \{N_j\}_{j=1}^{n} \in \mathbb{R}^n : \left\{ 0 < N_j^l \leq N_j \leq N_j^u,\ N_j^l < N_j^u \right\} \ \forall j \in \{1,\ldots,n\} \right\} \tag{17}$$

and the index sets defined in (16). Then, the multispecies population system (14) is <u>SOS</u> $F^{RN}$ iff:

$$\left\{\begin{array}{l} \left\{ \left\{ 0 \geq 1 - \alpha_{ii} N_i^u - \sum_{\substack{j \in \mathrm{A}_j^+ \\ j \neq i}} \alpha_{ij} N_j^l - \sum_{\substack{j \in \mathrm{A}_j^- \\ j \neq i}} \alpha_{ij} N_j^u \right\} \ \forall i \in R^+, \left\{ 0 \leq 1 - \alpha_{ii} N_i^l - \sum_{\substack{j \in \mathrm{A}_j^+ \\ j \neq i}} \alpha_{ij} N_j^u - \sum_{\substack{j \in \mathrm{A}_j^- \\ j \neq i}} \alpha_{ij} N_j^l \right\} \ \forall i \in R^+ \right\} \\ \left\{ \left\{ 0 \leq 1 - \alpha_{ii} N_i^u - \sum_{\substack{j \in \mathrm{A}_j^+ \\ j \neq i}} \alpha_{ij} N_j^u - \sum_{\substack{j \in \mathrm{A}_j^- \\ j \neq i}} \alpha_{ij} N_j^l \right\} \ \forall i \in R^-, \left\{ 0 \geq 1 - \alpha_{ii} N_i^l - \sum_{\substack{j \in \mathrm{A}_j^+ \\ j \neq i}} \alpha_{ij} N_j^l - \sum_{\substack{j \in \mathrm{A}_j^- \\ j \neq i}} \alpha_{ij} N_j^u \right\} \ \forall i \in R^- \right\} \end{array}\right\} \tag{18}$$

<u>Theorem 3</u>

Consider the rectangular set $F^{RN}$ defined in (17), the index sets defined in (16), and the admissible control strategy set $U^R$, and control variable functions $\{\alpha_{ii}\}_{i=1}^{n}$:

$$\begin{array}{l} U^R \triangleq \left\{ \{\hat{\alpha}_{ii}\}_{i=1}^{n} \in R^n : \left\{ 0 < \alpha_{ii}^l \leq \hat{\alpha}_{ii} \leq \alpha_{ii}^u,\ \alpha_{ii}^l < \alpha_{ii}^u \right\}, \forall i \in \{1,\ldots,n\} \right\} \\ \left\{ \{\alpha_{ii}\}_{i=1}^{n} : [0,\infty) \to U^R \subset \mathbb{R}^n,\ \{\alpha_{ii}\}_{i=1}^{n} : t \to \{\alpha_{ii}(t)\}_{i=1}^{n} \in U^R \subset \mathbb{R}^n \right\} \end{array} \tag{19}$$

Then, the multispecies population system (14) is <u>SIZOS</u> $F^{RN}$, $U^R$ iff:

$$\left\{\begin{array}{l}\left\{0 \ge 1-\alpha_{ii}^{u}N_i^u - \sum_{\substack{j\in A_i^+\\ j\ne i}} \alpha_{ij}N_j^l - \sum_{\substack{j\in A_i^-\\ j\ne i}} \alpha_{ij}N_j^u\right\} \ \forall i \in R^+, \left\{0 \le 1-\alpha_{ii}^{l}N_i^l - \sum_{\substack{j\in A_i^+\\ j\ne i}} \alpha_{ij}N_j^u - \sum_{\substack{j\in A_i^-\\ j\ne i}} \alpha_{ij}N_j^l\right\} \ \forall i \in R^+ \\ \left\{0 \le 1-\alpha_{ii}^{l}N_i^u - \sum_{\substack{j\in A_i^+\\ j\ne i}} \alpha_{ij}N_j^u - \sum_{\substack{j\in A_i^-\\ j\ne i}} \alpha_{ij}N_j^l\right\} \ \forall i \in R^-, \left\{0 \ge 1-\alpha_{ii}^{u}N_i^l - \sum_{\substack{j\in A_i^+\\ j\ne i}} \alpha_{ij}N_j^l - \sum_{\substack{j\in A_i^-\\ j\ne i}} \alpha_{ij}N_j^u\right\} \ \forall i \in R^-\end{array}\right\} \tag{20}$$

**3-species GLV model sustainability and sustainizability® case studies**

**Sustainability**

The sustainability of the 3 species GLV model presented by (May and Leonard, 1975) over a rectangular set is considered. The considered 3-species GLV system model is:

$$\left\{\begin{array}{l}\dfrac{dN_1(t)}{dt} = N_1(t)\left[1-N_1(t)-\alpha N_2(t)-\beta N_3(t)\right] \\ \dfrac{dN_2(t)}{dt} = N_2(t)\left[1-\beta N_1(t)-N_2(t)-\alpha N_3(t)\right] \\ \dfrac{dN_3(t)}{dt} = N_3(t)\left[1-\alpha N_1(t)-\beta N_2(t)-N_3(t)\right]\end{array}\right\} \tag{21}$$

The model arises from the general GLV model (14), as follows:

$$\left\{\begin{array}{l} n=3,\ r_1=r_2=r_3=1>0;\ \alpha_{11}=\alpha_{22}=\alpha_{33}=1>0 \\ \alpha_{12}=\alpha_{23}=\alpha_{31}=\alpha\ge\varepsilon_1>0;\ \alpha_{21}=\alpha_{32}=\alpha_{13}=\beta\ge\varepsilon_1>0\end{array}\right\} \Rightarrow \left\{\begin{array}{l}\left\{\begin{array}{l}A_i^+ \triangleq \{1,2,3\} \\ A_i^- \triangleq \varnothing\end{array}\right\} i=1,2,3 \\ R^+ \triangleq \{1,2,3\},\ R^- \triangleq \varnothing\end{array}\right\}.$$

As shown by (May and Leonard, 1975), the model possesses the following eight equilibria:

$$\begin{bmatrix}N_{1,e}\\N_{2,e}\\N_{3,e}\end{bmatrix} = \begin{bmatrix}0\\0\\0\end{bmatrix} \vee \begin{bmatrix}1\\0\\0\end{bmatrix} \vee \begin{bmatrix}0\\1\\0\end{bmatrix} \vee \begin{bmatrix}0\\0\\1\end{bmatrix} \vee \begin{bmatrix}\frac{1-\alpha}{1-\alpha\beta}\\ \frac{1-\beta}{1-\alpha\beta}\\ 0\end{bmatrix} \vee \begin{bmatrix}\frac{1-\beta}{1-\alpha\beta}\\ 0\\ \frac{1-\alpha}{1-\alpha\beta}\end{bmatrix} \vee \begin{bmatrix}0\\ \frac{1-\alpha}{1-\alpha\beta}\\ \frac{1-\beta}{1-\alpha\beta}\end{bmatrix} \vee \begin{bmatrix}\frac{1}{1+\alpha+\beta}\\ \frac{1}{1+\alpha+\beta}\\ \frac{1}{1+\alpha+\beta}\end{bmatrix}.$$

Only the last one of these eight equilibria features three nonzero populations, and for $\{\alpha>0;\ \beta>0\}$ the necessary and sufficient condition for its local stability is that $\alpha+\beta<2$.

The considered rectangular set is:

$$F^{R3} \triangleq \left\{\{N_j\}_{j=1}^3 \in \mathbb{R}^3 : \left\{N^l \le N_j \le N^u,\ 0<\varepsilon_2 \le N^l < N^u\right\}\ \forall j=1,2,3\right\} \tag{22}$$

Then applying Theorem 2 to $F^{R3}$ yields the necessary and sufficient <u>SOS</u> $F^{R3}$ conditions:

$$\left\{0 \ge \frac{1-N^u-N^l}{N^u} \ge (\alpha+\beta-1) \ge \frac{1-N^u-N^l}{N^l},\ 0<\varepsilon_1\le\alpha,\ 0<\varepsilon_1\le\beta,\ 0<\varepsilon_2\le N^l<N^u\right\} \tag{23}$$

It is clear that (23) cannot possibly be satisfied if $\alpha+\beta>1$ or $N^u+N^l<1$ or $N^l>1$.

Considering the model's competition coefficients $(\alpha,\beta)\in[\varepsilon_1,\infty)\times[\varepsilon_1,\infty)$ to be known and fixed and $\varepsilon_1$ to be infinitesimally small, the above <u>SOS</u> $F^{R3}$ conditions (23) then identify the maximal region of population bounds $\left(N^l, N^u\right) \in \mathbb{R}^2$, such that system (21) is SOS $F^{R3}$, as the hatched triangle, not including its left boundary, that is illustrated in Figure 1a.
Considering the rectangular $F^{R3}$ set's population bounds $\left(N^l, N^u\right) \in \mathbb{R}^2$ to be known and fixed, the SOS $F^{R3}$ conditions (23) identify the maximal region of the model's competition coefficients $(\alpha, \beta) \in [\varepsilon_1, \infty) \times [\varepsilon_1, \infty)$ such that system (21) is SOS $F^{R3}$, as the hatched trapezoid in Figure 1b, which also demonstrates that $\frac{1-N^l}{N^u} \geq \frac{1-N^u}{N^l}$ must hold for the maximal region to be nonempty.

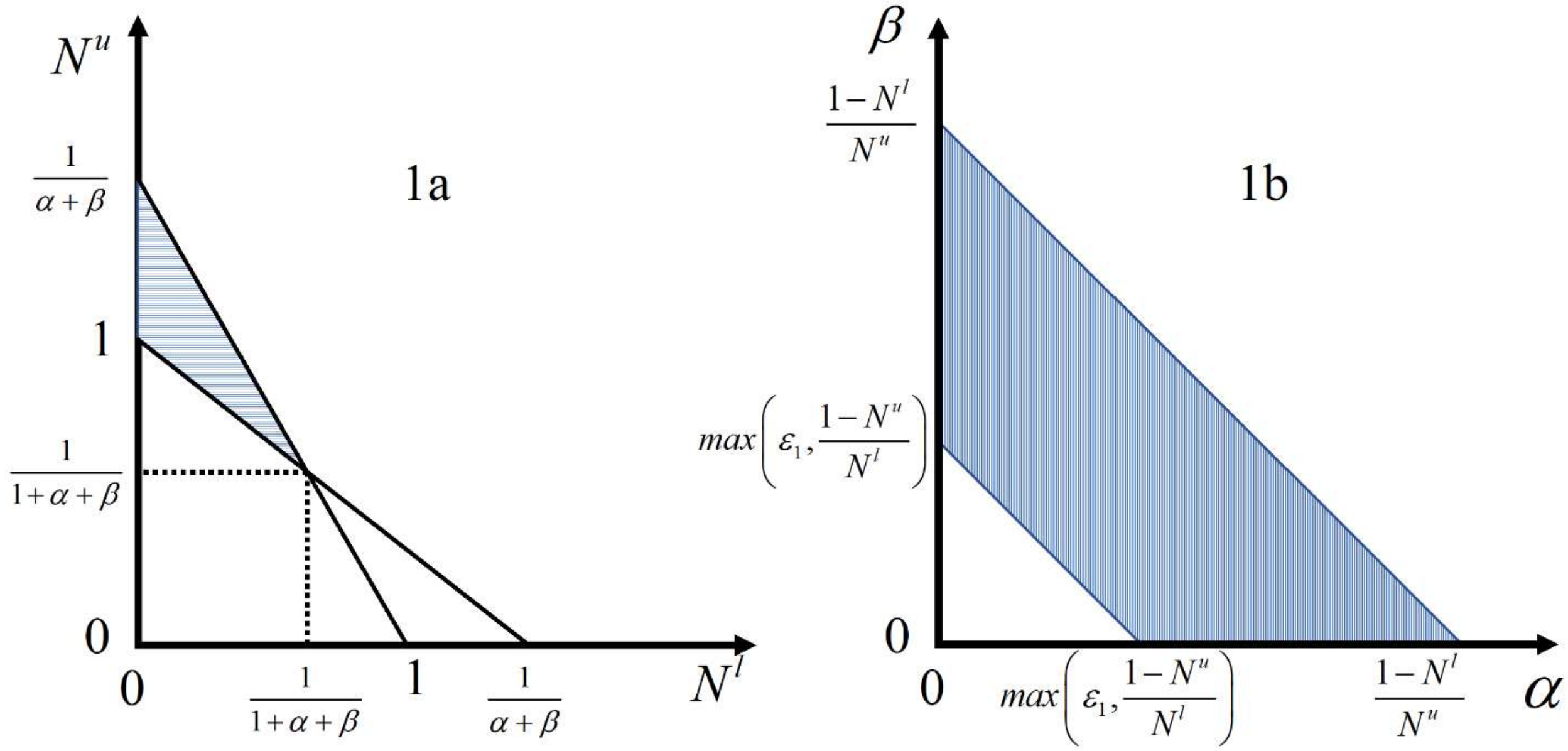

Figure 1. Maximal Population Bound $\left(N^l, N^u\right) \in \mathbb{R}^2$ (1a) and Maximal Competition Coefficient $(\alpha, \beta) \in [\varepsilon_1, \infty) \times [\varepsilon_1, \infty)$ (1b) Regions for (21) to be SOS $F^{R3}$

The above SOS $F^{R3}$ conditions (23) are not satisfied, if and only if the system's vector field points outward in some portions of the $F^{R3}$ set's boundary, i.e. if and only if

$$\exists(N_1,N_2,N_3):\left\{\begin{array}{l}\left\{N_2 \in \left(N^l,N^u\right) \wedge N_3 \in \left(N^l,N^u\right) \wedge \left\{\begin{array}{l}\left\{N_1 = N^l \wedge 1-N^l-\alpha N_2-\beta N_3<0\right\} \vee \\ \left\{N_1 = N^u \wedge 1-N^u-\alpha N_2-\beta N_3>0\right\}\end{array}\right\}\right\} \vee \\ \left\{N_1 \in \left(N^l,N^u\right) \wedge N_3 \in \left(N^l,N^u\right) \wedge \left\{\begin{array}{l}\left\{N_2 = N^l \wedge 1-\beta N_1-N^l-\alpha N_3<0\right\} \vee \\ \left\{N_2 = N^u \wedge 1-\beta N_1-N^u-\alpha N_3>0\right\}\end{array}\right\}\right\} \vee \\ \left\{N_1 \in \left(N^l,N^u\right) \wedge N_2 \in \left(N^l,N^u\right);; \wedge \left\{\begin{array}{l}\left\{N_3 = N^l \wedge 1-\alpha N_1-\beta N_2-N^l<0\right\} \vee \\ \left\{N_3 = N^u \wedge 1-\alpha N_1-\beta N_2-N^u>0\right\}\end{array}\right\}\right\}\end{array}\right\} \qquad (24)$$

To demonstrate the validity and value of these theoretical results derived based on Theorem 2, next we consider the following numerical case studies:
Case study 1
Consider model (21) with fixed parameter values $\alpha = 0.2,\ \beta = 0.05$. Then the maximal triangle of population bounds $\left(N^l, N^u\right) \in \mathbb{R}^2$ so that (21) is SOS $F^{R3}$, as illustrated in Figure 1, has vertices

$$\left(N^l,N^u\right)=(0,1),\left(N^l,N^u\right)=\left(0,\frac{1}{(\alpha+\beta)}\right)=(0,4),\left(N^l,N^u\right)=\left(\frac{1}{(1+\alpha+\beta)},\frac{1}{(1+\alpha+\beta)}\right)=(0.8,0.8)$$

Case study 1a

Consider the rectangular set $F^{R3}$ with $\left(N^l,N^u\right)=(0.5,2.0)$, which is a point that belongs to the right boundary of the maximal sustainable population bound triangular region illustrated in Figure 1. Further, none of the conditions in equation (24) are satisfied, and thus the system's vector field either points inward or is tangential to the $F^{R3}$ set's boundary. This implies that (21) is SOS $F^{R3}$. As illustrated in Figure 2, in the state-space $\left\{N_j\right\}_{j=1}^{3}\in\mathbb{R}^3$ the system (21) with parameter values $\alpha=0.2,\ \beta=0.05$ features the stable equilibrium point $\left(N_1,N_2,N_3\right)=\left((1+\alpha+\beta)^{-1},(1+\alpha+\beta)^{-1},(1+\alpha+\beta)^{-1}\right)=(0.8,0.8,0.8)$, since $\alpha+\beta=0.25<2$, which is shown to attract all trajectories initiated at the eight vertices of the set $F^{R3}$.

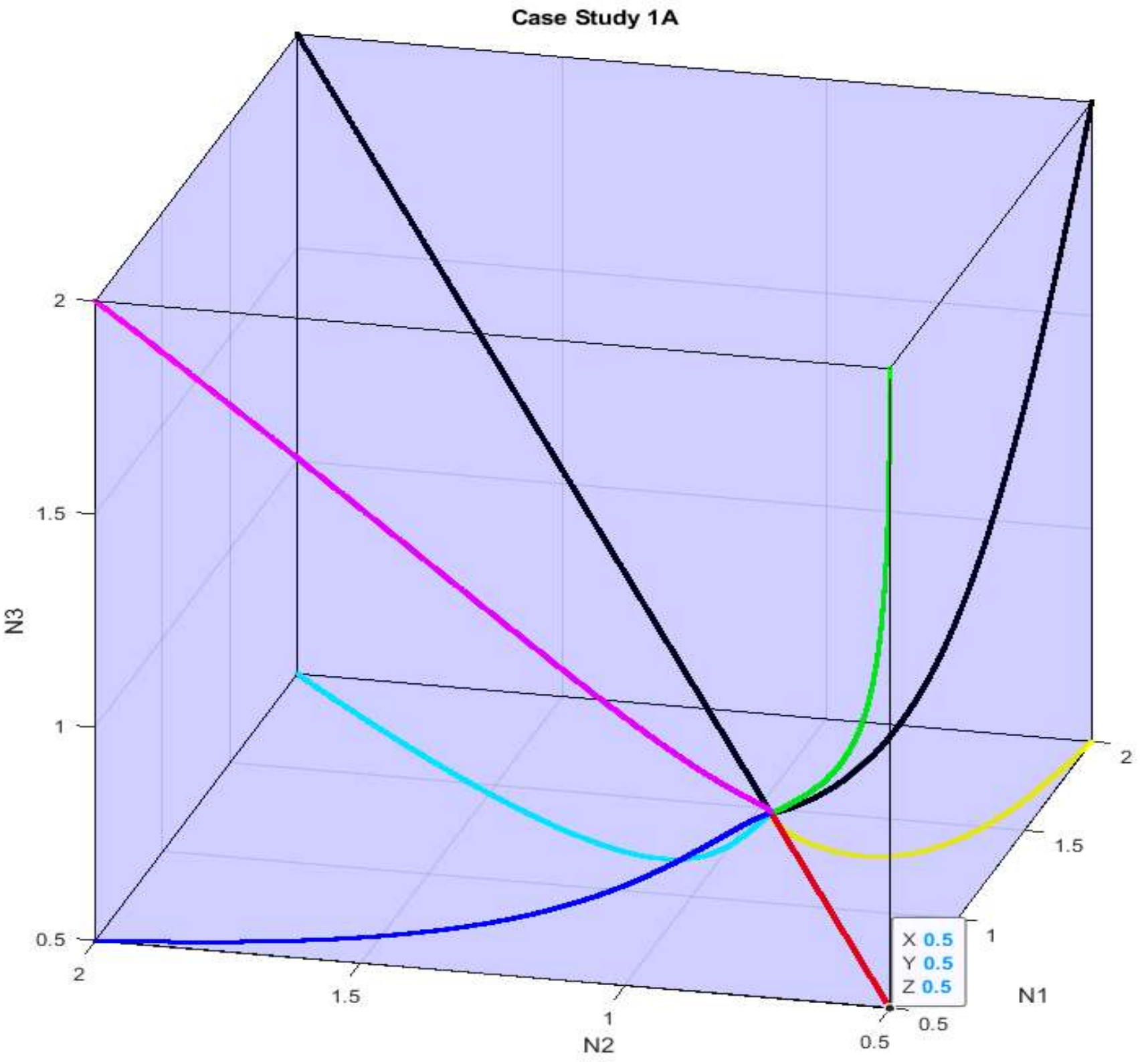

Figure 2. System (21), $(\alpha,\beta)=(0.2,0.05);\left(N^l,N^u\right)=(0.5,2.0)$ is SOS $F^{R3}$

All trajectories initiated at the eight $F^{R3}$ vertices are attracted to the stable equilibrium point $\left(N_1,N_2,N_3\right)=(0.8,0.8,0.8)$, and remain forever within $F^{R3}$

The time evolution of the 3 species involved in the trajectories initiated at the two $F^{R3}$ vertices $(N_1, N_2, N_3) = (N^l, N^l, N^u) = (0.5, 0.5, 2.0)$, $(N_1, N_2, N_3) = (N^l, N^u, N^l) = (0.5, 2.0, 0.5)$ is captured in Figure 3, which illustrates convergence to the equilibrium point and that the populations remain for all time within $F^{R3}$.

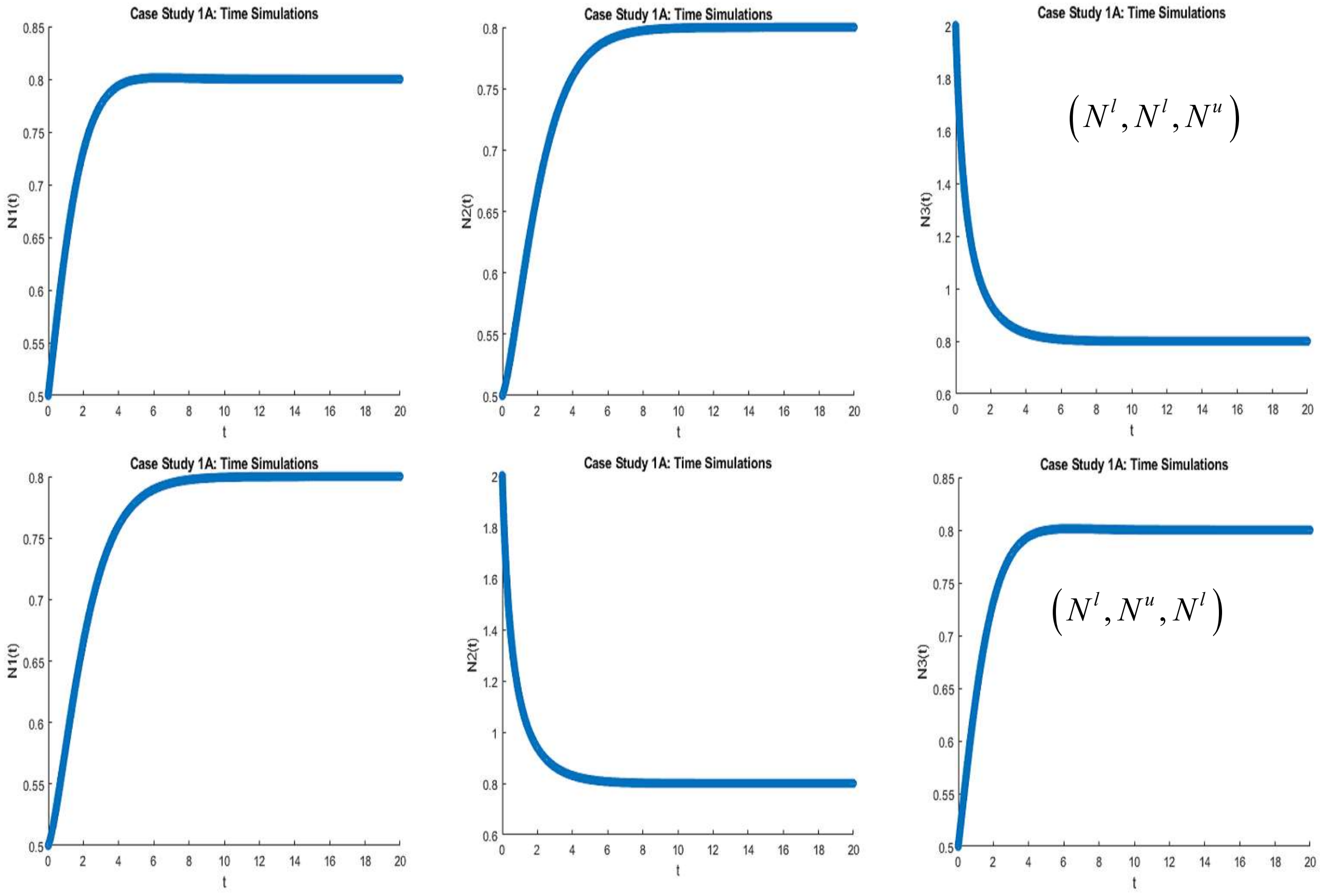

Figure 3. Time evolution of system (21), $(\alpha, \beta) = (0.2, 0.05)$, $(N^l, N^u) = (0.5, 2.0)$ 3 species trajectories initiated at $(N_1, N_2, N_3) = (N^l, N^l, N^u) = (0.5, 0.5, 2.0)$ (top row), $(N_1, N_2, N_3) = (N^l, N^u, N^l) = (0.5, 2.0, 0.5)$ (bottom row)

<u>Case study 1b</u>

Consider the rectangular set $F^{R3}$ with $(N^l, N^u) = (0.75, 3.25)$, which is a point that lies outside the maximal sustainable population bound triangular region illustrated in Figure 1. Further, some of the conditions in equation (24) are satisfied, and thus the system's vector field points outward in some portions of the $F^{R3}$ set's boundary. This implies that (21) is not <u>SOS</u> $F^{R3}$. This is confirmed in Figure 4, which illustrates in the state-space $\{N_j\}_{j=1}^{3} \in \mathbb{R}^3$ that although the system (21) with parameter values $\alpha = 0.2,\ \beta = 0.05$ still features the stable equilibrium point $(N_1, N_2, N_3) = \left((1+\alpha+\beta)^{-1}, (1+\alpha+\beta)^{-1}, (1+\alpha+\beta)^{-1}\right) = (0.8, 0.8, 0.8)$, since $\alpha + \beta = 0.25 < 2$, some of the trajectories initiated at the $F^{R3}$ set's eight vertices do not remain for all time within $F^{R3}$, even though they are eventually attracted to the aforementioned stable equilibrium point.

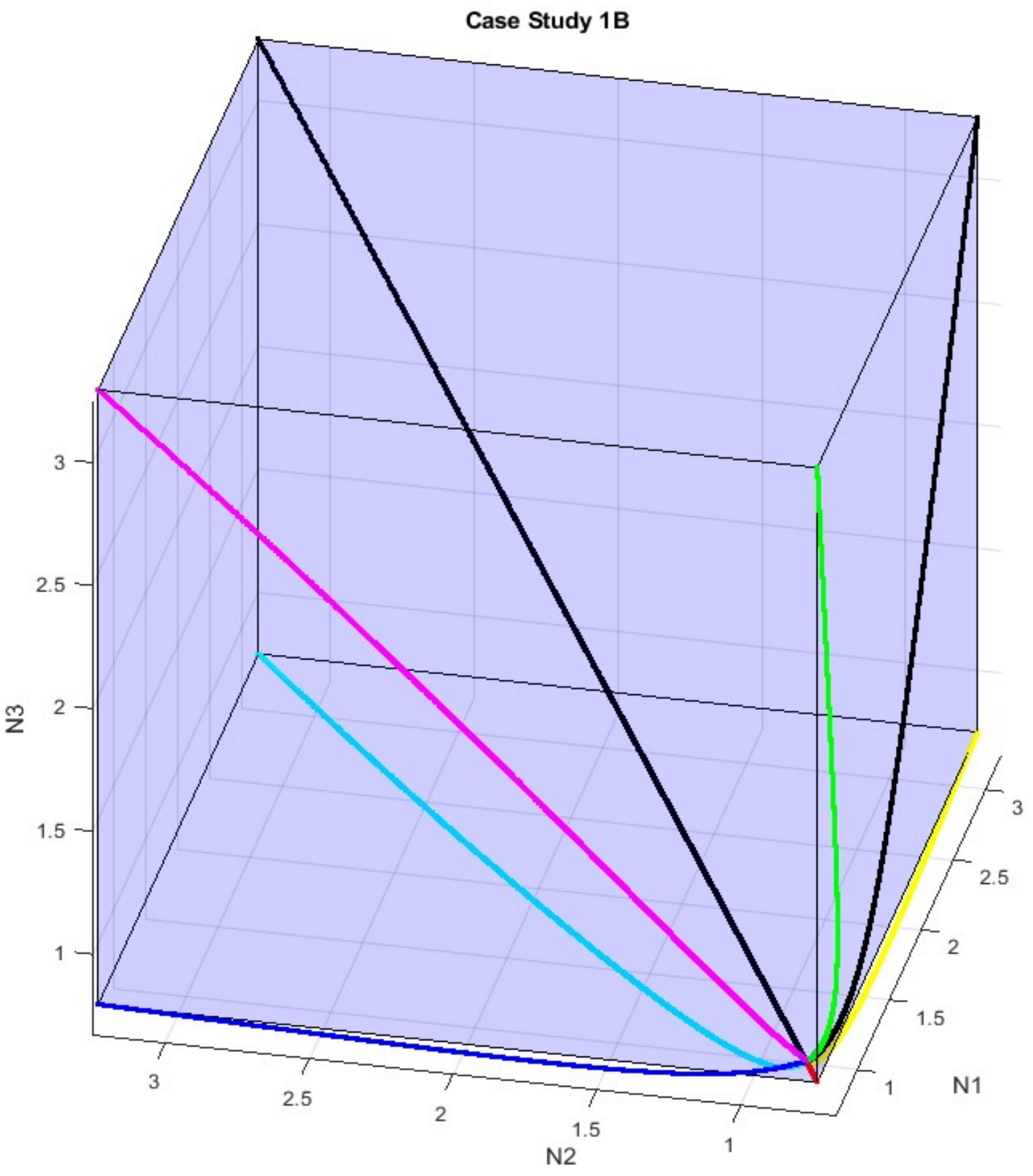

Figure 4. System (21), $(\alpha,\beta)=(0.2,0.05);(N^l,N^u)=(0.75,3.25)$ isn't SOS $F^{R3}$

All trajectories initiated at the eight $F^{R3}$ vertices are attracted to the stable equilibrium point $(N_1,N_2,N_3)=(0.8,0.8,0.8)$, but some do not remain forever within $F^{R3}$

The time evolution of the 3 species involved in the trajectories initiated at the two $F^{R3}$ vertices $(N_1,N_2,N_3)=(N^l,N^l,N^u)=(0.75,0.75,3.25)$, $(N_1,N_2,N_3)=(N^l,N^u,N^l)=(0.75,3.25,0.75)$,

is captured in Figure 5, which illustrates convergence to the equilibrium point, but the populations do not remain forever within $F^{R3}$.

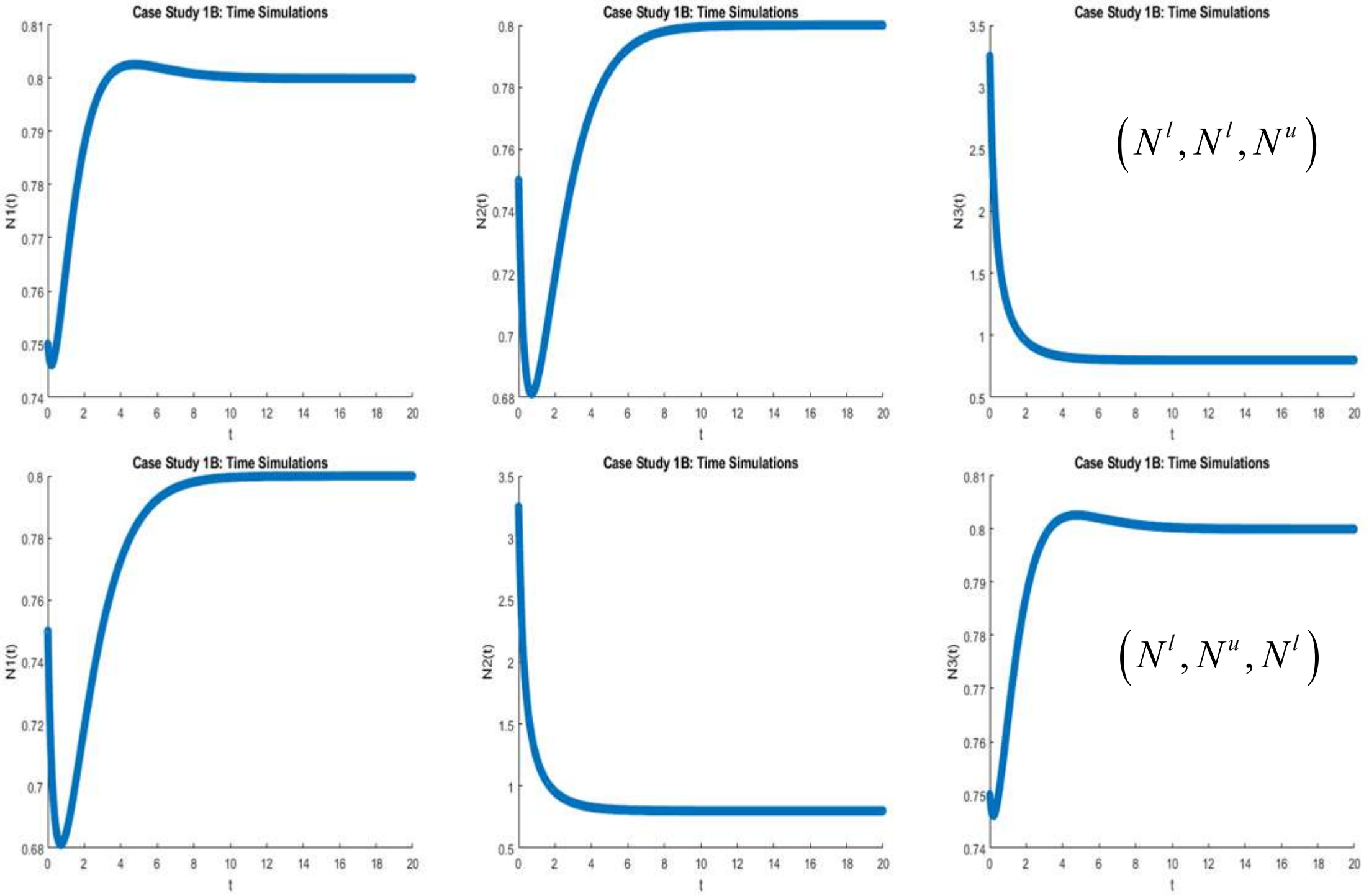

Figure 5. Time evolution of system (21), $(\alpha, \beta) = (0.2, 0.05), (N^l, N^u) = (0.75, 3.25)$ 3 species trajectories initiated at $(N_1, N_2, N_3) = (N^l, N^l, N^u) = (0.75, 0.75, 3.25)$ (top row), $(N_1, N_2, N_3) = (N^l, N^u, N^l) = (0.75, 3.25, 0.75)$ (bottom row)

<u>Case study 2</u>

Consider model (21) with fixed parameter values $\alpha = 0.8,\ \beta = 1.3$. As proved earlier, the <u>SOS</u> $F^{R3}$ conditions (23), illustrated in Figure 1, cannot possibly be satisfied since $\alpha + \beta > 1$. Further, some of the conditions in equation (24) are satisfied, and thus the system's vector field points outward in some portions of the $F^{R3}$ set's boundary. This implies that system (21) cannot be <u>SOS</u> $F^{R3}$ for any set $F^{R3}$ of the form defined in (22). This is confirmed in Figure 6 for the set $F^{R3}$ with $(N^l, N^u) = (0.25, 0.38)$, which illustrates in the state-space $\{N_j\}_{j=1}^{3} \in \mathbb{R}^3$ that all trajectories initiated at the eight vertices of $F^{R3}$, are not attracted by its unstable equilibrium point, ($\alpha + \beta = 2.1 > 2$),

$$(N_1, N_2, N_3) = \left( \frac{1}{1+\alpha+\beta}, \frac{1}{1+\alpha+\beta}, \frac{1}{1+\alpha+\beta} \right) = \left( \frac{1}{3.1}, \frac{1}{3.1}, \frac{1}{3.1} \right) = (0.32258, 0.32258, 0.32258),$$

exhibit chaotic behavior, and do not remain within $F^{R3}$.

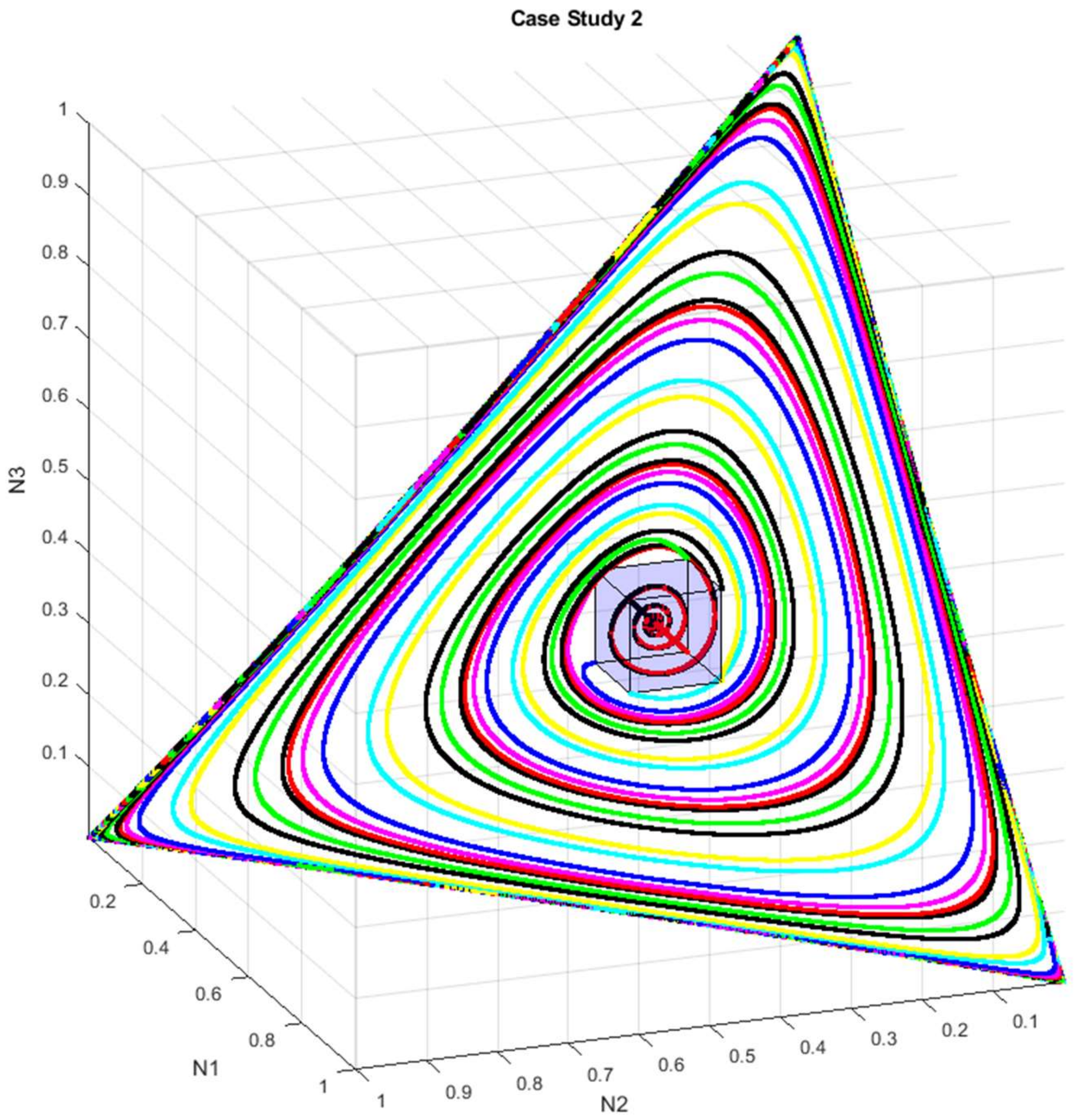

Figure 6. System (21), $(\alpha,\beta)=(0.8,1.3);(N^l,N^u)=(0.25,0.38)$ isn't SOS $F^{R3}$

All trajectories initiated at the eight $F^{R3}$ vertices are not attracted to the unstable equilibrium point $(N_1,N_2,N_3)=(1/3.1,1/3.1,1/3.1)$, exhibit chaotic behavior and do not remain within $F^{R3}$ The time evolution of the 3 species involved in the trajectories initiated at the two $F^{R3}$ vertices $(N_1,N_2,N_3)=(N^l,N^l,N^l)=(0.25,0.25,0.25)$, $(N_1,N_2,N_3)=(N^l,N^l,N^u)=(0.25,0.25,0.38)$, is captured in Figure 7, which confirms the aforementioned behavior.

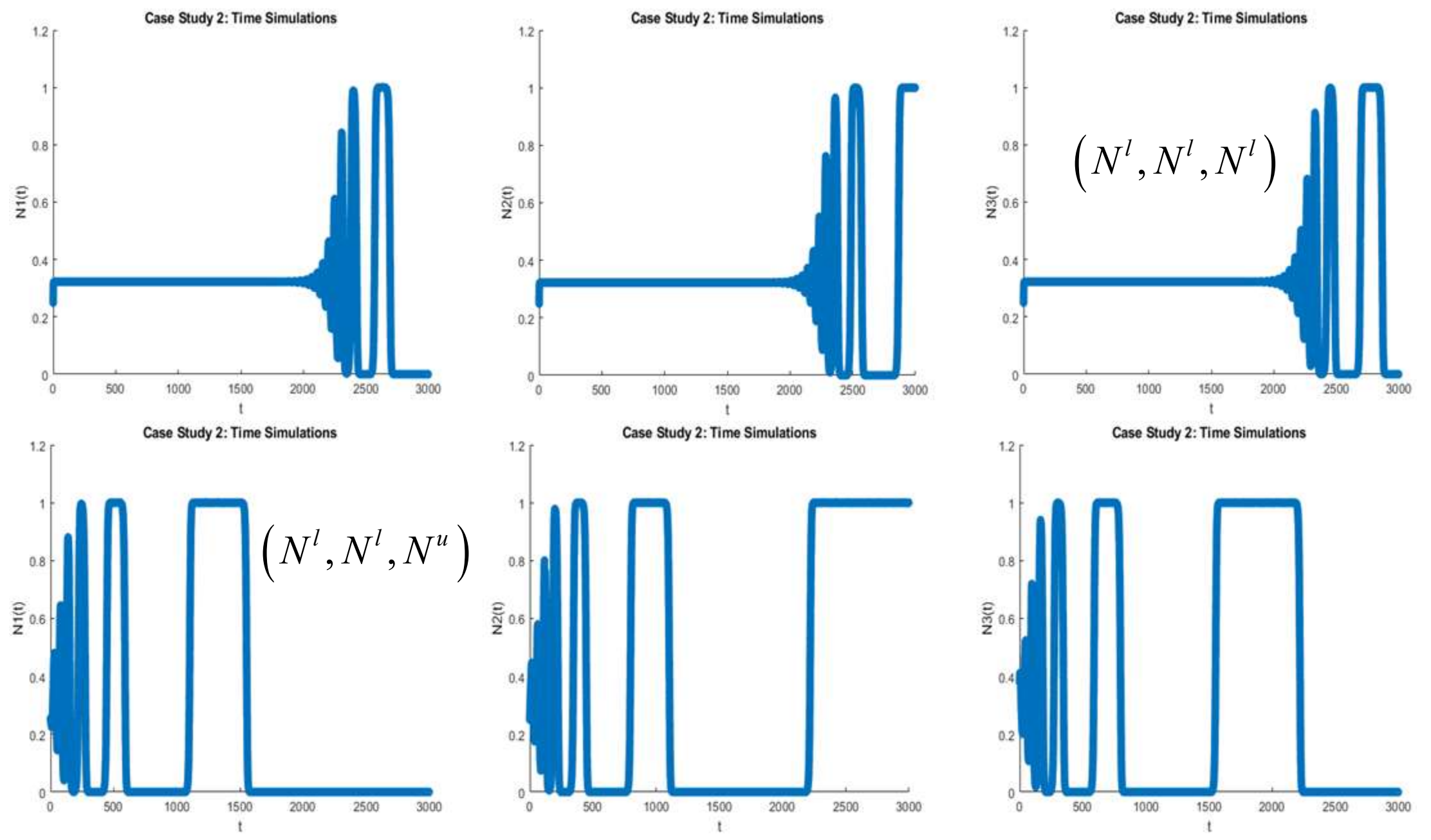

Figure 7. Time evolution of System (21), $(\alpha,\beta)=(0.8,1.3),\left(N^l,N^u\right)=(0.25,0.38)$ 3 species trajectories initiated at $\left(N_1,N_2,N_3\right)=\left(N^l,N^l,N^l\right)=(0.25,0.25,0.25)$ (top row), $\left(N_1,N_2,N_3\right)=\left(N^l,N^l,N^u\right)=(0.25,0.25,0.38)$ (bottom row)

**Sustainizability®**

The sustainizability® of the 3 species GLV model presented by (May and Leonard, 1975) over a rectangular set is considered.

The considered 3-species GLV model is:

$$\left\{\begin{aligned} \frac{dN_1(t)}{dt} &= N_1(t)\left[1-\alpha_{11}(t)N_1(t)-\alpha N_2(t)-\beta N_3(t)\right] \\ \frac{dN_2(t)}{dt} &= N_2(t)\left[1-\beta N_1(t)-\alpha_{22}(t)N_2(t)-\alpha N_3(t)\right] \\ \frac{dN_3(t)}{dt} &= N_3(t)\left[1-\alpha N_1(t)-\beta N_2(t)-\alpha_{33}(t)N_3(t)\right] \end{aligned}\right\} \tag{25}$$

which arises from the general GLV model (14), as follows:

$$\left\{\begin{array}{l} n=3,\ r_1=r_2=r_3=1>0; \\ \alpha_{12}=\alpha_{23}=\alpha_{31}=\alpha\geq\varepsilon_1>0;\ \alpha_{21}=\alpha_{32}=\alpha_{13}=\beta\geq\varepsilon_1>0 \end{array}\right\} \Rightarrow \left\{\begin{array}{l} \left\{\begin{array}{l} \mathrm{A}_i^+ \triangleq \{1,2,3\} \\ \mathrm{A}_i^- \triangleq \varnothing \end{array}\right\} i=1,2,3 \\ R^+ \triangleq \{1,2,3\},\ R^- \triangleq \varnothing \end{array}\right\},$$

Consider the GLV model (14) as a forced system, with the self-competition coefficients being employed as control variable functions $\left\{\alpha_{ii}\right\}_{i=1}^3 : \left\{[0,\infty)\rightarrow U^R\subset\mathbb{R}^3,\ t\rightarrow\left\{\alpha_{ii}(t)\right\}_{i=1}^3\in U^R\subset\mathbb{R}^3\right\}$, with admissible control strategy set

$$U^R \triangleq \left\{\left\{\hat{\alpha}_{ii}\right\}_{i=1}^3\in\mathbb{R}^3 : \left\{0<\alpha^l\leq\hat{\alpha}_{ii}\leq\alpha^u\ \ \forall i=1,n\right\}\right\} \tag{26}$$

Then, the above definition of the set $U^R$, the definition of the set $F^{R3}$ in (22), and the SIZOS conditions established in Theorem 3, yield the following necessary and sufficient conditions for the above 3-species GLV model (25) to be SIZOS $F^{R3}$, $U^R$:

$$\left\{0 \ge 1-\alpha^u N^u - (\alpha+\beta)N^l,\ 0 \le 1-\alpha^l N^l - (\alpha+\beta)N^u,\ 0<\varepsilon_1 \le \alpha,\ 0<\varepsilon_1 \le \beta,\ 0<\varepsilon_2 \le N^l < N^u, \alpha^l \le \alpha^u\right\} \Leftrightarrow$$

$$\left\{\begin{bmatrix} \dfrac{1-\alpha^l N^l}{N^u} \ge (\alpha+\beta) \ge \dfrac{1-\alpha^u N^u}{N^l} \\ \alpha^u \ge \dfrac{N^u - N^l\left(1-\alpha^l N^l\right)}{\left(N^u\right)^2} \end{bmatrix}, \alpha^u \ge \frac{1-(\alpha+\beta)N^l}{N^u}, \alpha^l \le \frac{1-(\alpha+\beta)N^u}{N^l}, \begin{bmatrix} 0<\varepsilon_1 \le \alpha \\ 0<\varepsilon_1 \le \beta \\ 0<\varepsilon_2 \le N^l < N^u \\ \alpha^l \le \alpha^u \end{bmatrix}\right\} \tag{27}$$

If $\alpha^l = \alpha^u = 1$,then there is effectively no control strategy and the 3-species GLV model (25) dynamic behavior is identical to the dynamic behavior of the 3-species GLV model (21) discussed in the Sustainability case studies 1, 2, presented above.

In general however $\alpha^l \le \alpha^u$, and then admissible control strategies can be applied that can alter the 3-species GLV model (25) dynamic behavior from that of the 3-species GLV model (21). It is then clear that the above conditions cannot possibly be satisfied if $\alpha^l N^l > 1$, which in turn implies that it must hold $\alpha^u \ge \alpha^l > \frac{1}{N^l} > 0$. To demonstrate the validity and value of these theoretical results derived based on Theorem 3, next we consider the following case study:

Case study 3

Consider model (25) with the fixed parameter values $\alpha = 0.8,\ \beta = 1.3$ used in Case study 2 for the unforced system (21). Then, as discussed in Case study 2, the unforced system (21) cannot be SOS $F^R$ for any set $F^R$ of the form defined in (19), and its three nonzero population equilibrium point $(N_1, N_2, N_3) = (0.32258, 0.32258, 0.32258)$, is not stable, since $\alpha+\beta = 2.1 > 2$.

However, consider model (25) with the fixed parameter values $\alpha = 0.8,\ \beta = 1.3$, and the rectangular set $F^{R3}$ with $(N^l, N^u) = (0.25, 0.38)$, and the admissible control strategy set $U^R \triangleq \left\{\left\{\hat{\alpha}_{ii}\right\}_{i=1}^3 \in R^3 : \left\{0 < \alpha^l \le \hat{\alpha}_{ii} \le \alpha^u\ \ \forall i = 1, n\right\}\right\}$. Then, the above derived SIZOS $F^{R3}$, $U^R$ conditions (27) are satisfied if and only if :

$$\left\{\alpha^u \ge \frac{1-(\alpha+\beta)N^l}{N^u} = \frac{1-2.1\cdot 0.25}{0.38} = 1.25, \alpha^l \le \frac{1-(\alpha+\beta)N^u}{N^l} = \frac{1-2.1\cdot 0.38}{0.25} = 0.808, \alpha^l \le \alpha^u\right\}.$$

This suggests that there exists a control strategy within the admissible control strategy set $U^R \triangleq \left\{\left\{\hat{\alpha}_{ii}\right\}_{i=1}^3 \in R^3 : \left\{0 < 0.808 = \alpha^l \le \hat{\alpha}_{ii} \le \alpha^u = 1.25\ \forall i = 1, n\right\}\right\}$ with $(\alpha^l, \alpha^u) = (0.808, 1.25)$ that will make the forced system (25) SOS $F^{R3}$. A particular such control strategy that lies within $U^R$ is:

$$\bar{\alpha}_{ii} : \left\{N_j\right\}_{j=1}^3 \to \bar{\alpha}_{ii}\left(\left\{N_j\right\}_{j=1}^3\right) \triangleq \left[0.808 + 192\cdot\begin{bmatrix} max(0, N_i - 0.250) - \\ max(0, N_i - 0.251) \end{bmatrix} + 250\cdot\begin{bmatrix} max(0, N_i - 0.379) - \\ max(0, N_i - 0.380) \end{bmatrix}\right] =$$

$$= \begin{bmatrix} 0.808\ if\ N_i \le 0.250 \\ 0.808 + 192\cdot(N_i - 0.250) \le 1\ if\ N_i \in (0.250, 0.251] \\ 1\ if\ N_i \in (0.251, 0.379] \\ 1 + 250\cdot(N_i - 0.379) \le 1.25\ if\ N_i \in (0.379, 0.380] \\ 1.25\ if\ N_i > 0.380 \end{bmatrix} \in \left[\alpha^l, \alpha^u\right] = [0.808, 1.25]\ \forall i \in \{1,2,3\} \tag{28}$$

The unforced system that results from the use of this control strategy is:

$$\left\{\begin{aligned}
\frac{dN_1(t)}{dt} &= N_1(t)\left[1-\left[\begin{array}{l}0.808+192\cdot\begin{bmatrix}max\left(0,N_1(t)-0.250\right)-\\ max\left(0,N_1(t)-0.251\right)\end{bmatrix}\\ +250\cdot\begin{bmatrix}max\left(0,N_1(t)-0.379\right)-\\ max\left(0,N_1(t)-0.380\right)\end{bmatrix}\end{array}\right]N_1(t)-0.8N_2(t)-1.3N_3(t)\right]\\
\frac{dN_2(t)}{dt} &= N_2(t)\left[1-1.3N_1(t)-\left[\begin{array}{l}0.808+192\cdot\begin{bmatrix}max\left(0,N_2(t)-0.250\right)-\\ max\left(0,N_2(t)-0.251\right)\end{bmatrix}\\ +250\cdot\begin{bmatrix}max\left(0,N_2(t)-0.379\right)-\\ max\left(0,N_2(t)-0.380\right)\end{bmatrix}\end{array}\right]N_2(t)-0.8N_3(t)\right]\\
\frac{dN_3(t)}{dt} &= N_3(t)\left[1-0.8N_1(t)-1.3N_2(t)-\left[\begin{array}{l}0.808+192\cdot\begin{bmatrix}max\left(0,N_3(t)-0.250\right)-\\ max\left(0,N_3(t)-0.251\right)\end{bmatrix}\\ +250\cdot\begin{bmatrix}max\left(0,N_3(t)-0.379\right)-\\ max\left(0,N_3(t)-0.380\right)\end{bmatrix}\end{array}\right]N_3(t)\right]
\end{aligned}\right\} \tag{29}$$

The unforced system (29) also features the same unstable equilibrium point $\left(N_1,N_2,N_3\right)=\left(0.32258,0.32258,0.32258\right)$, as the system (21) considered in case study 2 with fixed parameter values $\alpha=0.8,\ \beta=1.3$ .However, unlike the system (21) in case study 2, which was not SOS $F^{R3}$, the system (29) is SOS $F^{R3}$, and its vector field either points inward or is tangential to the $F^{R3}$ set's boundary, which implies that the 3-species GLV model (25) is SIZOS $F^{R3}$, $U^R$, confirming the validity of Theorem 2. This is confirmed in Figure 8, which illustrates in the state-space $\left\{N_j\right\}_{j=1}^{3}\in\mathbb{R}^3$ that the trajectories of the unforced system (29) initiated at the eight vertices of the rectangular set $F^{R3}$ are not attracted to the unstable equilibrium point $\left(N_1,N_2,N_3\right)=\left(0.32258,0.32258,0.32258\right)$, remain within $F^{R3}$, although they remain chaotic. This confirms the 3-species GLV model (25) is SIZOS $F^{R3}$, $U^R$, and that the SIZOS methodology can deliver controlled chaos.

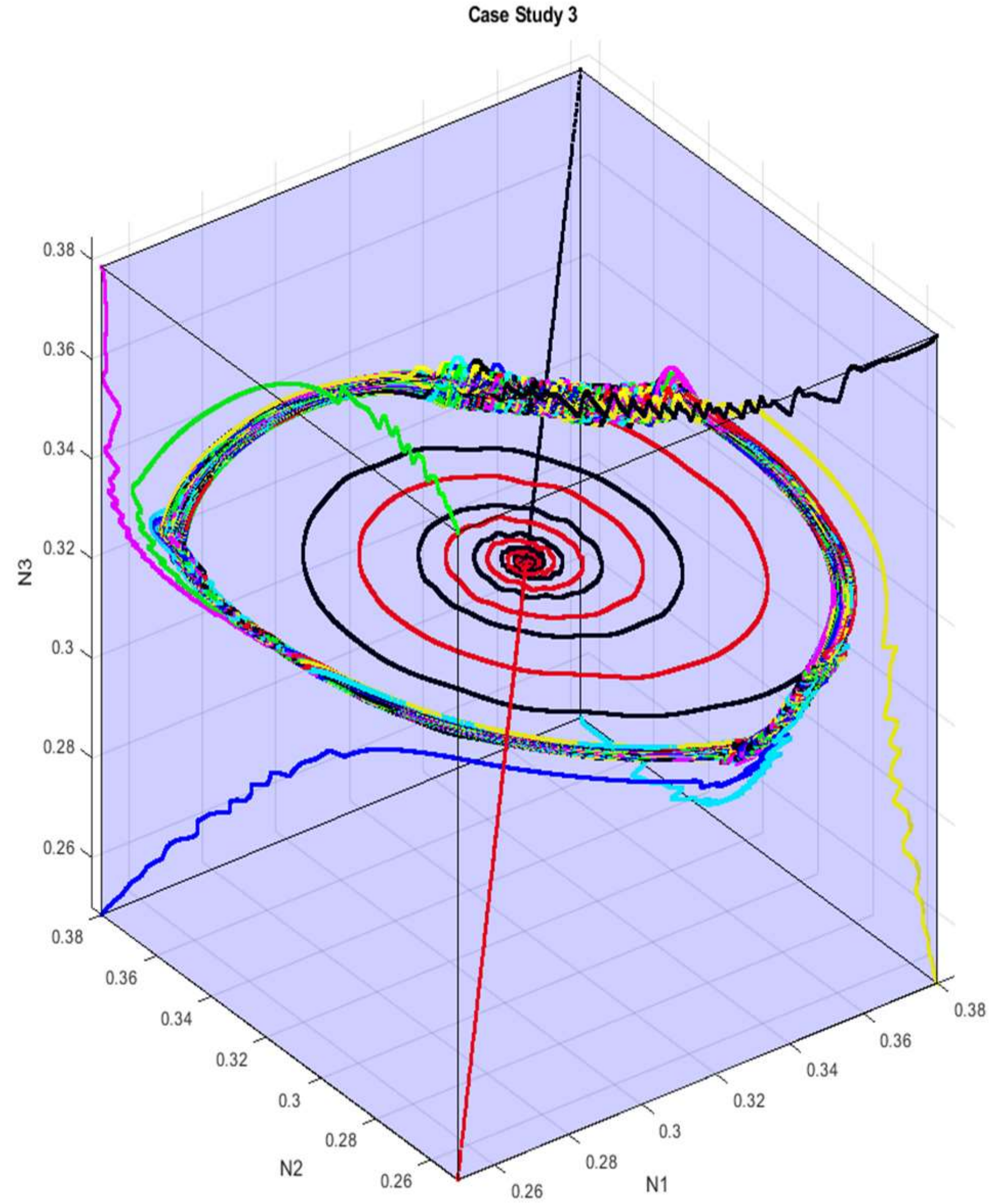

Figure 8. System (29), $(\alpha,\beta)=(0.8,1.3);(N^l,N^u)=(0.25,0.38)$ is <u>SOS</u> $F^{R3}$ which confirms the 3-species GLV system (25) is <u>SIZOS</u> $F^{R3}$, $U^R$

All trajectories of (31) initiated at the eight $F^{R3}$ vertices are not attracted to the unstable equilibrium point $(N_1,N_2,N_3)=(1/3.1,1/3.1,1/3.1)$, are chaotic, but remain forever within $F^{R3}$

The time evolution of the 3 species involved in the trajectories initiated at the $F^{R3}$ vertex $(N_1,N_2,N_3)=(N^l,N^l,N^l)=(0.25,0.25,0.25)$ is captured in Figure 9, which illustrates no convergence to the equilibrium point, chaotic behavior, and all species remaining within $F^{R3}$.

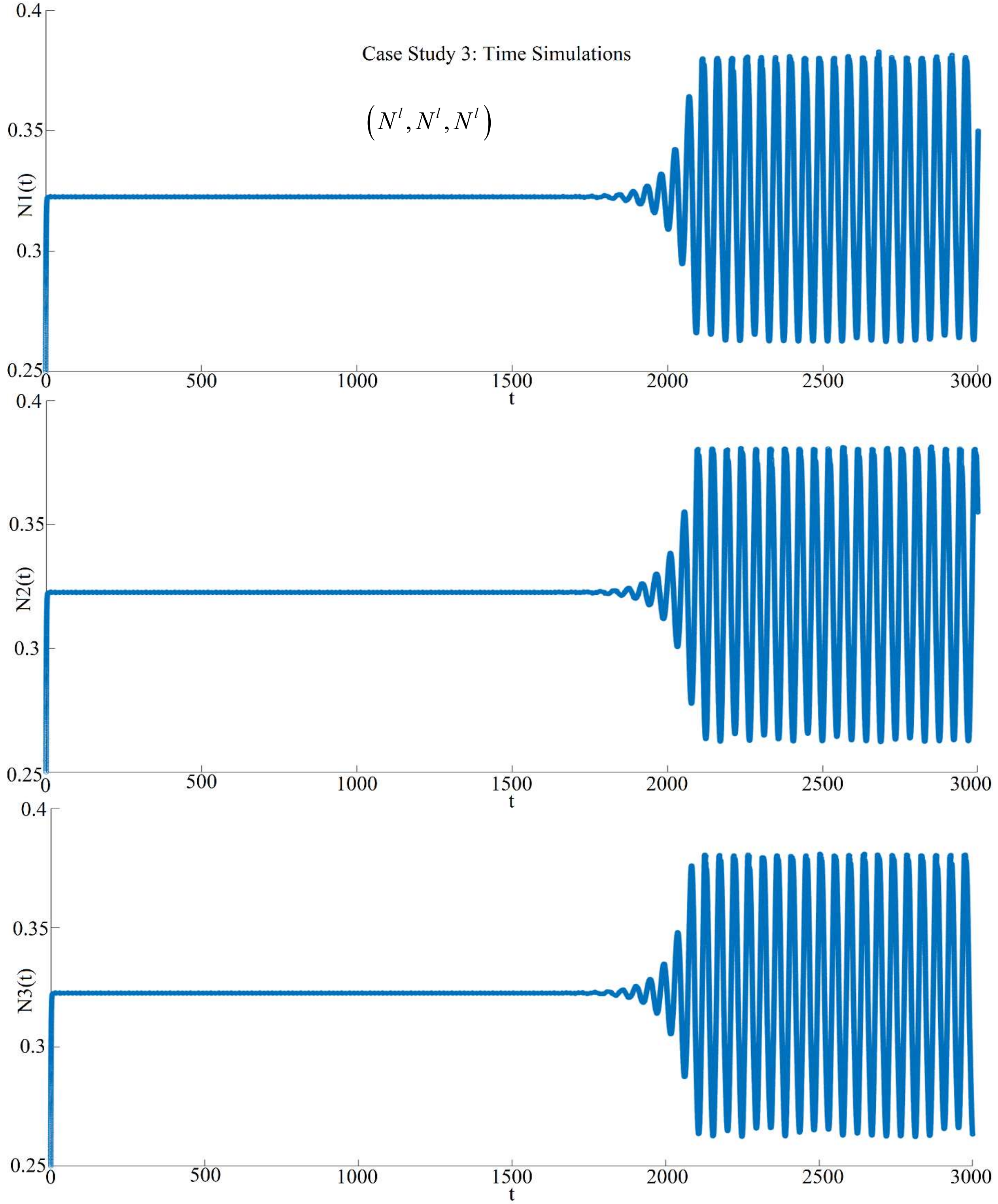

Figure 9. Time evolution of System (29), $(\alpha, \beta) = (0.8, 1.3), \left(N^l, N^u\right) = (0.25, 0.38)$ 3 species trajectories initiated at $\left(N_1, N_2, N_3\right) = \left(N^l, N^l, N^l\right) = (0.25, 0.25, 0.25)$, not attracted to equilibrium point $\left(N_1, N_2, N_3\right) = (1/3.1, 1/3.1, 1/3.1)$, exhibit chaos but stay within $F^{R3}$

## Conclusions

A Theorem has been proved that establishes a game theory based necessary and sufficient sustainizability® condition for a general ODE system. Subsequently, two other theorems are proved for the general, n-species Gause-Lotka-Volterra (GLV) population model, establishing necessary and sufficient sustainability and sustainizability® conditions over rectangular sets. Three case studies on the May-Leonard, 3-species, GLV model are then presented that illustrate the power of the aforementioned Theorems. It was demonstrated that the Theorems' predictions are validated through numerical simulations, and even chaotic systems with unsustainable behavior can be made sustainable through control strategies identified by the proposed methodologies.

**Title: A Game Theoretic Approach to Sustainizability® Over Sets and its application to a multi-species population model**
**Order of Authors:** Ioannis V. Manousiouthakis[a], Vasilios I. Manousiouthakis[b*]
**First Author:** Ioannis V. Manousiouthakis[a], [a] Hydrogen Engineering Research Company, LLC (H-E-R-C, LLC), Los Angeles, CA, 90077, United States of America
**Corresponding Author:** Vasilios I. Manousiouthakis[b], [b] Department of Chemical and Biomolecular Engineering, Hydrogen Engineering Research Consortium (HERC), University of California at Los Angeles (UCLA), Los Angeles, CA, USA (ORCID: 0000-0002-5926-9923)
Email: vasilios@ucla.edu, Phone: (310) 206-0300, Address: 5549 Boelter Hall, Box 951592, Los Angeles, CA, 90095-1592, U.S.A.

**Supplementary Information**

Proof of Theorem 1

Application of the SIZOS $F,U$ conditions of equation (10) yields that the ODE system of Theorem 1 is SIZOS $F,U$ iff:

$$\left\{\forall \{z_j\}_{j=1}^n : \left\{\begin{array}{l} \{z_j\}_{j=1}^n \in F \\ S\left(\{z_j\}_{j=1}^n\right) \neq \emptyset \end{array}\right\} \exists \{\bar{u}_l\}_{l=1}^p : \left\{\begin{array}{l} \{\bar{u}_l\}_{l=1}^p : F \to U,\ \{\bar{u}_l\}_{l=1}^p : \{x_j\}_{j=1}^n \to \left\{\bar{u}_l\left(\{x_j\}_{j=1}^n\right)\right\}_{l=1}^p \in U \\ \sum_{i=1}^n \left[\frac{\partial \Phi_k\left(\{z_j\}_{j=1}^n\right)}{\partial z_i} \cdot f_i\left(\{z_j\}_{j=1}^n, \left\{\bar{u}_l\left(\{z_j\}_{j=1}^n\right)\right\}_{l=1}^p\right)\right] \le 0\ \forall k \in S\left(\{z_j\}_{j=1}^n\right) \end{array}\right\}\right\}.$$

Consider an arbitrary but fixed $\{z_j\}_{j=1}^n \in F$ such that $S\left(\{z_j\}_{j=1}^n\right) \neq \emptyset$.

It then holds:
$$\left\{\exists \{\bar{u}_l\}_{l=1}^p : \left\{\begin{array}{l} \{\bar{u}_l\}_{l=1}^p : F \to U,\ \{\bar{u}_l\}_{l=1}^p : \{z_j\}_{j=1}^n \to \left\{\bar{u}_l\left(\{z_j\}_{j=1}^n\right)\right\}_{l=1}^p \in U \\ \sum_{i=1}^n \left[\frac{\partial \Phi_k\left(\{z_j\}_{j=1}^n\right)}{\partial z_i} \cdot f_i\left(\{z_j\}_{j=1}^n, \left\{\bar{u}_l\left(\{z_j\}_{j=1}^n\right)\right\}_{l=1}^p\right)\right] \le 0\ \forall k \in S\left(\{z_j\}_{j=1}^n\right) \neq \emptyset \end{array}\right\}\right\} \Leftrightarrow$$

$$\left\{\exists \{\bar{\bar{u}}_l\}_{l=1}^p : \left\{\{\bar{\bar{u}}_l\}_{l=1}^p \in U,\ \sum_{i=1}^n \left[\frac{\partial \Phi_k\left(\{z_j\}_{j=1}^n\right)}{\partial z_i} \cdot f_i\left(\{z_j\}_{j=1}^n, \{\bar{\bar{u}}_l\}_{l=1}^p\right)\right] \le 0\ \forall k \in S\left(\{z_j\}_{j=1}^n\right) \neq \emptyset\right\}\right\} \Leftrightarrow$$

$$0 \ge \min_{\{\bar{\bar{u}}_l\}_{l=1}^p \in U} \max_{k \in S\left(\{z_j\}_{j=1}^n\right)} \sum_{i=1}^n \left[\frac{\partial \Phi_k\left(\{z_j\}_{j=1}^n\right)}{\partial z_i} \cdot f_i\left(\{z_j\}_{j=1}^n, \{\bar{\bar{u}}_l\}_{l=1}^p\right)\right].$$

Then the ODE system of Theorem 1 is SIZOS $F,U$ iff:

$$\left\{0 \geq \min_{\{\bar{\bar{u}}_l\}_{l=1}^p \in U} \max_{k \in S\left(\{z_j\}_{j=1}^n\right)} \sum_{i=1}^n \left[\frac{\partial \Phi_k\left(\{z_j\}_{j=1}^n\right)}{\partial z_i} \cdot f_i\left(\{z_j\}_{j=1}^n, \{\bar{\bar{u}}_l\}_{l=1}^p\right)\right] \forall \{z_j\}_{j=1}^n : \begin{Bmatrix} \{z_j\}_{j=1}^n \in F \\ S\left(\{z_j\}_{j=1}^n\right) \neq \varnothing \end{Bmatrix}\right\} \Leftrightarrow$$

$$\left\{0 \geq \max_{\substack{\{z_j\}_{j=1}^n \in F \\ S\left(\{z_j\}_{j=1}^n\right) \neq \varnothing}} \min_{\{\bar{\bar{u}}_l\}_{l=1}^p \in U} \max_{k \in S\left(\{z_j\}_{j=1}^n\right)} \sum_{i=1}^n \left[\frac{\partial \Phi_k\left(\{z_j\}_{j=1}^n\right)}{\partial z_i} \cdot f_i\left(\{z_j\}_{j=1}^n, \{\bar{\bar{u}}_l\}_{l=1}^p\right)\right]\right\} \text{O.E.}\Delta.$$

Proof of Theorem 2

Application of the SOS $F^R$ conditions of equation (6) to this Theorem, yields that the multispecies population ODE system (14) is SOS $F^{RN}$ iff:

$$\left\{\begin{array}{l} \left\{r_i N_i^{i,u}\left[1-\sum_{j=1}^n \alpha_{ij} N_j^{i,u}\right] \leq 0 \;\forall \{N_j^{i,u}\}_{j=1}^n : \begin{bmatrix} N_i^{i,u} = N_i^u \\ N_j^l \leq N_j^{i,u} \leq N_j^u \;\; \forall j \in \{1,\ldots,n\}; j \neq i \end{bmatrix}\right\} \;\forall i \in \{1,\ldots,n\} \\ \left\{r_i N_i^{i,l}\left[1-\sum_{j=1}^n \alpha_{ij} N_j^{i,l}\right] \geq 0 \;\forall \{N_j^{i,l}\}_{j=1}^n : \begin{bmatrix} N_i^{i,l} = N_i^l \\ N_j^l \leq N_j^{i,l} \leq N_j^u \;\; \forall j \in \{1,\ldots,n\}; j \neq i \end{bmatrix}\right\} \;\forall i \in \{1,\ldots,n\} \end{array}\right\} \underset{N_i^l > 0 \;\forall i \in \{1,\ldots,n\}}{\Leftrightarrow}$$

$$\left\{\begin{array}{l} \left\{1-\sum_{j=1}^n \alpha_{ij} N_j^{i,u} \leq 0 \;\forall \{N_j^{i,u}\}_{j=1}^n : \begin{bmatrix} N_i^{i,u} = N_i^u \\ N_j^l \leq N_j^{i,u} \leq N_j^u \;\; \forall j \in \{1,\ldots,n\}; j \neq i \end{bmatrix}\right\} \;\forall i \in R^+ \\ \left\{1-\sum_{j=1}^n \alpha_{ij} N_j^{i,u} \geq 0 \;\forall \{N_j^{i,u}\}_{j=1}^n : \begin{bmatrix} N_i^{i,u} = N_i^u \\ N_j^l \leq N_j^{i,u} \leq N_j^u \;\; \forall j \in \{1,\ldots,n\}; j \neq i \end{bmatrix}\right\} \;\forall \in R^- \\ \left\{1-\sum_{j=1}^n \alpha_{ij} N_j^{i,l} \geq 0 \;\forall \{N_j^{i,l}\}_{j=1}^n : \begin{bmatrix} N_i^{i,l} = N_i^l \\ N_j^l \leq N_j^{i,l} \leq N_j^u \;\; \forall j \in \{1,\ldots,n\}; j \neq i \end{bmatrix}\right\} \;\forall i \in R^+ \\ \left\{1-\sum_{j=1}^n \alpha_{ij} N_j^{i,l} \leq 0 \;\forall \{N_j^{i,l}\}_{j=1}^n : \begin{bmatrix} N_i^{i,l} = N_i^l \\ N_j^l \leq N_j^{i,l} \leq N_j^u \;\; \forall j \in \{1,\ldots,n\}; j \neq i \end{bmatrix}\right\} \;\forall i \in R^- \end{array}\right\} \Leftrightarrow$$

$$\left\{\begin{array}{l} \left\{0 \geq \left\{\begin{array}{l} \max_{\{N_j^{i,u}\}_{j=1}^n} \left[1-\sum_{j=1}^n \alpha_{ij} N_j^{i,u}\right] \\ s.t.\;\; N_i^{i,u} = N_i^u,\; N_j^l \leq N_j^{i,u} \leq N_j^u \;\; \forall j \in \{1,\ldots,n\}; j \neq i \end{array}\right\}\right\} \;\forall i \in R^+ \\ \left\{0 \leq \left\{\begin{array}{l} \min_{\{N_j^{i,u}\}_{j=1}^n} \left[1-\sum_{j=1}^n \alpha_{ij} N_j^{i,u}\right] \\ s.t.\;\; N_i^{i,u} = N_i^u,\; N_j^l \leq N_j^{i,u} \leq N_j^u \;\; \forall j \in \{1,\ldots,n\}; j \neq i \end{array}\right\}\right\} \;\forall i \in R^- \\ \left\{0 \leq \left\{\begin{array}{l} \min_{\{N_j^{i,l}\}_{j=1}^n} \left[1-\sum_{j=1}^n \alpha_{ij} N_j^{i,l}\right] \\ s.t.\;\; N_i^{i,l} = N_i^l,\; N_j^l \leq N_j^{i,l} \leq N_j^u \;\; \forall j \in \{1,\ldots,n\}; j \neq i \end{array}\right\}\right\} \;\forall i \in R^+ \\ \left\{0 \geq \left\{\begin{array}{l} \max_{\{N_j^{i,l}\}_{j=1}^n} \left[1-\sum_{j=1}^n \alpha_{ij} N_j^{i,l}\right] \\ s.t.\;\; N_i^{i,l} = N_i^l,\; N_j^l \leq N_j^{i,l} \leq N_j^u \;\; \forall j \in \{1,\ldots,n\}; j \neq i \end{array}\right\}\right\} \;\forall i \in R^- \end{array}\right\} \Leftrightarrow$$

$$\left\{\begin{array}{l}\left\{0\geq 1-\alpha_{ii}N_i^u-\sum_{\substack{j\in A_j^+\\ j\neq i}}\alpha_{ij}N_j^l-\sum_{\substack{j\in A_j^-\\ j\neq i}}\alpha_{ij}N_j^u\right\}\ \forall i\in R^+\\ \left\{0\leq 1-\alpha_{ii}N_i^u-\sum_{\substack{j\in A_j^+\\ j\neq i}}\alpha_{ij}N_j^u-\sum_{\substack{j\in A_j^-\\ j\neq i}}\alpha_{ij}N_j^l\right\}\ \forall i\in R^-\\ \left\{0\leq 1-\alpha_{ii}N_i^l-\sum_{\substack{j\in A_j^+\\ j\neq i}}\alpha_{ij}N_j^u-\sum_{\substack{j\in A_j^-\\ j\neq i}}\alpha_{ij}N_j^l\right\}\ \forall i\in R^+\\ \left\{0\geq 1-\alpha_{ii}N_i^l-\sum_{\substack{j\in A_j^+\\ j\neq i}}\alpha_{ij}N_j^l-\sum_{\substack{j\in A_j^-\\ j\neq i}}\alpha_{ij}N_j^u\right\}\ \forall i\in R^-\end{array}\right\}\text{ O.E.}\Delta.$$

Proof of Theorem 3

Application of the SIZOS $F^R$, $U^R$ conditions of equations (10), (13) to this Theorem, yields that the multispecies population ODE system (14) is SIZOS $F^{RN}$, $U^R$ iff:

$$\left\{\left\{\forall\{N_j\}_{j=1}^n:\left[\begin{array}{c}\left\{\begin{array}{c}S^U\left(\{N_j\}_{j=1}^n\right)\neq\varnothing\\ \vee\\ S^L\left(\{N_j\}_{j=1}^n\right)\neq\varnothing\end{array}\right\}\\ \{N_j\}_{j=1}^n\in\left\{\left[N_j^l,N_j^u\right]\right\}_{i=1}^n\end{array}\right]\right\}\exists\{\bar{\alpha}_{ii}\}_{i=1}^n:\left\{\begin{array}{l}\{\bar{\alpha}_{ii}\}_{i=1}^n:\{N_j\}_{j=1}^n\to\left\{\bar{\alpha}_{ii}\left(\{N_j\}_{j=1}^n\right)\right\}_{i=1}^n\in\left\{\left[\alpha_{ii}^l,\alpha_{ii}^u\right]\right\}_{i=1}^n\\ \forall i\in S^U\left(\{N_j\}_{j=1}^n\right)\ r_iN_i\left[\begin{array}{l}1-\bar{\alpha}_{ii}\left(\{N_j\}_{j=1}^n\right)N_i-\\ -\sum_{\substack{j=1\\ j\neq i}}^n\alpha_{ij}N_j\end{array}\right]\leq 0\\ \forall i\in S^L\left(\{N_j\}_{j=1}^n\right)\ r_iN_i\left[\begin{array}{l}1-\bar{\alpha}_{ii}\left(\{N_j\}_{j=1}^n\right)N_i-\\ -\sum_{\substack{j=1\\ j\neq i}}^n\alpha_{ij}N_j\end{array}\right]\geq 0\end{array}\right\}\right\}\underset{\forall i\in\{1,\dots,n\}}{\overset{N_i^l>0}{\Leftrightarrow}}$$

$$\left\{\left\{\forall\{N_j\}_{j=1}^n:\left[\begin{array}{c}\left\{\begin{array}{c}S^U\left(\{N_j\}_{j=1}^n\right)\neq\varnothing\\ \vee\\ S^L\left(\{N_j\}_{j=1}^n\right)\neq\varnothing\end{array}\right\}\\ \{N_j\}_{j=1}^n\in\left\{\left[N_j^l,N_j^u\right]\right\}_{i=1}^n\end{array}\right]\right\}\exists\{\bar{\alpha}_{ii}\}_{i=1}^n:\left\{\begin{array}{l}\{\bar{\alpha}_{ii}\}_{i=1}^n:\{N_j\}_{j=1}^n\to\left\{\bar{\alpha}_{ii}\left(\{N_j\}_{j=1}^n\right)\right\}_{i=1}^n\in\left\{\left[\alpha_{ii}^l,\alpha_{ii}^u\right]\right\}_{i=1}^n\\ \forall i\in S^U\left(\{N_j\}_{j=1}^n\right)\ r_i\left[1-\bar{\alpha}_{ii}\left(\{N_j\}_{j=1}^n\right)N_i-\sum_{\substack{j=1\\ j\neq i}}^n\alpha_{ij}N_j\right]\leq 0\\ \forall i\in S^L\left(\{N_j\}_{j=1}^n\right)\ r_i\left[1-\bar{\alpha}_{ii}\left(\{N_j\}_{j=1}^n\right)N_i-\sum_{\substack{j=1\\ j\neq i}}^n\alpha_{ij}N_j\right]\geq 0\end{array}\right\}\right\}\Leftrightarrow$$

$$\left\{\left\{\forall\left\{N_j\right\}_{j=1}^n:\begin{bmatrix}\begin{Bmatrix} S^U\left(\left\{N_j\right\}_{j=1}^n\right)\neq\varnothing\\ \vee \\ S^L\left(\left\{N_j\right\}_{j=1}^n\right)\neq\varnothing\end{Bmatrix}\\ \left\{N_j\right\}_{j=1}^n\in\left\{\left[N_j^l,N_j^u\right]\right\}_{i=1}^n\end{bmatrix}\right\}\exists\left\{\bar{\alpha}_{ii}\right\}_{i=1}^n:\left\{\begin{array}{l} \left\{\bar{\alpha}_{ii}\right\}_{i=1}^n:\left\{N_j\right\}_{j=1}^n\to\left\{\bar{\alpha}_{ii}\left(\left\{N_j\right\}_{j=1}^n\right)\right\}_{i=1}^n\in\left\{\left[\alpha_{ii}^l,\alpha_{ii}^u\right]\right\}_{i=1}^n \\ \forall i\in S^U\left(\left\{N_j\right\}_{j=1}^n\right)\cap R^+ \;\; 1-\bar{\alpha}_{ii}\left(\left\{N_j\right\}_{j=1}^n\right)N_i-\sum_{\substack{j=1\\ j\neq i}}^n \alpha_{ij}N_j\le 0 \\ \forall i\in S^U\left(\left\{N_j\right\}_{j=1}^n\right)\cap R^- \;\; 1-\bar{\alpha}_{ii}\left(\left\{N_j\right\}_{j=1}^n\right)N_i-\sum_{\substack{j=1\\ j\neq i}}^n \alpha_{ij}N_j\ge 0 \\ \forall i\in S^L\left(\left\{N_j\right\}_{j=1}^n\right)\cap R^+ \;\; 1-\bar{\alpha}_{ii}\left(\left\{N_j\right\}_{j=1}^n\right)N_i-\sum_{\substack{j=1\\ j\neq i}}^n \alpha_{ij}N_j\ge 0 \\ \forall i\in S^L\left(\left\{N_j\right\}_{j=1}^n\right)\cap R^- \;\; 1-\bar{\alpha}_{ii}\left(\left\{N_j\right\}_{j=1}^n\right)N_i-\sum_{\substack{j=1\\ j\neq i}}^n \alpha_{ij}N_j\le 0 \end{array}\right\}\right\}$$

Consider an arbitrary $\left\{N_j\right\}_{j=1}^n\in\left\{\left[N_j^l,N_j^u\right]\right\}_{i=1}^n:S^U\left(\left\{N_j\right\}_{j=1}^n\right)\cup S^L\left(\left\{N_j\right\}_{j=1}^n\right)\neq\varnothing$, and the corresponding index sets $S^U\left(\left\{N_j\right\}_{j=1}^n\right)\cap R^+, S^U\left(\left\{N_j\right\}_{j=1}^n\right)\cap R^-, S^L\left(\left\{N_j\right\}_{j=1}^n\right)\cap R^+, S^L\left(\left\{N_j\right\}_{j=1}^n\right)\cap R^-$ which are mutually exclusive given the definition (17) of the rectangular set $F^{RN}$, the index set definitions (16), and that $N_j^l<N_j^u \;\; \forall j=1,\ldots,n$, $R^+\cap R^-=\varnothing$. Further, to each $i$ that is an element of the union of these index sets, corresponds one and only one function $\bar{\alpha}_{ii}$ which has the same index subscript, and appears only in the ith right hand side of equation (14), and only in one of the four maximized bracketed expressions above. This implies that the selection of this function $\bar{\alpha}_{ii}$ will not influence any of the aforementioned bracketed expressions, except the one in which the $\bar{\alpha}_{ii}$ appears. Therefore, the above condition is equivalent to:

$$\left\{\left\{\forall\left\{N_j\right\}_{j=1}^n:\begin{bmatrix}\begin{Bmatrix} S^U\left(\left\{N_j\right\}_{j=1}^n\right)\neq\varnothing\\ \vee \\ S^L\left(\left\{N_j\right\}_{j=1}^n\right)\neq\varnothing\end{Bmatrix}\\ \left\{N_j\right\}_{j=1}^n\in\left\{\left[N_j^l,N_j^u\right]\right\}_{i=1}^n\end{bmatrix}\right\}\left\{\begin{array}{l} \forall i\in S^U\left(\left\{N_j\right\}_{j=1}^n\right)\cap R^+ \;\exists\bar{\alpha}_{ii}:\left\{\begin{array}{l}\left\{N_j\right\}_{j=1}^n\to\bar{\alpha}_{ii}\left(\left\{N_j\right\}_{j=1}^n\right)\in\left[\alpha_{ii}^l,\alpha_{ii}^u\right]\\ 1-\bar{\alpha}_{ii}\left(\left\{N_j\right\}_{j=1}^n\right)N_i-\sum_{\substack{j=1\\ j\neq i}}^n \alpha_{ij}N_j\le 0\end{array}\right\} \\ \forall i\in S^U\left(\left\{N_j\right\}_{j=1}^n\right)\cap R^- \;\exists\bar{\alpha}_{ii}:\left\{\begin{array}{l}\left\{N_j\right\}_{j=1}^n\to\bar{\alpha}_{ii}\left(\left\{N_j\right\}_{j=1}^n\right)\in\left[\alpha_{ii}^l,\alpha_{ii}^u\right]\\ 1-\bar{\alpha}_{ii}\left(\left\{N_j\right\}_{j=1}^n\right)N_i-\sum_{\substack{j=1\\ j\neq i}}^n \alpha_{ij}N_j\ge 0\end{array}\right\} \\ \forall i\in S^L\left(\left\{N_j\right\}_{j=1}^n\right)\cap R^+ \;\exists\bar{\alpha}_{ii}:\left\{\begin{array}{l}\left\{N_j\right\}_{j=1}^n\to\bar{\alpha}_{ii}\left(\left\{N_j\right\}_{j=1}^n\right)\in\left[\alpha_{ii}^l,\alpha_{ii}^u\right]\\ 1-\bar{\alpha}_{ii}\left(\left\{N_j\right\}_{j=1}^n\right)N_i-\sum_{\substack{j=1\\ j\neq i}}^n \alpha_{ij}N_j\ge 0\end{array}\right\} \\ \forall i\in S^L\left(\left\{N_j\right\}_{j=1}^n\right)\cap R^- \;\exists\bar{\alpha}_{ii}:\left\{\begin{array}{l}\left\{N_j\right\}_{j=1}^n\to\bar{\alpha}_{ii}\left(\left\{N_j\right\}_{j=1}^n\right)\in\left[\alpha_{ii}^l,\alpha_{ii}^u\right]\\ 1-\bar{\alpha}_{ii}\left(\left\{N_j\right\}_{j=1}^n\right)N_i-\sum_{\substack{j=1\\ j\neq i}}^n \alpha_{ij}N_j\le 0\end{array}\right\} \end{array}\right\}\right\}\Leftrightarrow$$

$$
\left\{\begin{array}{l}
\left\{\begin{array}{l}\forall i\in R^{+}\\ \forall\left\{N_j\right\}_{j=1}^{n}:\begin{cases}N_i=N_i^{u}\\ N_j\in\left[N_j^{l},N_j^{u}\right]\ \forall j\in\left\{1,n\right\}-\left\{i\right\}\end{cases}\end{array}\right\}\exists\bar{\alpha}_{ii}:\left\{\begin{array}{l}\left\{N_j\right\}_{j=1}^{n}\to\bar{\alpha}_{ii}\left(\left\{N_j\right\}_{j=1}^{n}\right)\in\left[\alpha_{ii}^{l},\alpha_{ii}^{u}\right]\\ 1-\bar{\alpha}_{ii}\left(\left\{N_j\right\}_{j=1}^{n}\right)N_i^{u}-\sum_{\substack{j=1\\ j\neq i}}^{n}\alpha_{ij}N_j\le 0\end{array}\right\}\\
\left\{\begin{array}{l}\forall i\in R^{-}\\ \forall\left\{N_j\right\}_{j=1}^{n}:\begin{cases}N_i=N_i^{u}\\ N_j\in\left[N_j^{l},N_j^{u}\right]\ \forall j\in\left\{1,n\right\}-\left\{i\right\}\end{cases}\end{array}\right\}\exists\bar{\alpha}_{ii}:\left\{\begin{array}{l}\left\{N_j\right\}_{j=1}^{n}\to\bar{\alpha}_{ii}\left(\left\{N_j\right\}_{j=1}^{n}\right)\in\left[\alpha_{ii}^{l},\alpha_{ii}^{u}\right]\\ 1-\bar{\alpha}_{ii}\left(\left\{N_j\right\}_{j=1}^{n}\right)N_i^{u}-\sum_{\substack{j=1\\ j\neq i}}^{n}\alpha_{ij}N_j\ge 0\end{array}\right\}\\
\left\{\begin{array}{l}\forall i\in R^{+}\\ \forall\left\{N_j\right\}_{j=1}^{n}:\begin{cases}N_i=N_i^{l}\\ N_j\in\left[N_j^{l},N_j^{u}\right]\ \forall j\in\left\{1,n\right\}-\left\{i\right\}\end{cases}\end{array}\right\}\exists\bar{\alpha}_{ii}:\left\{\begin{array}{l}\left\{N_j\right\}_{j=1}^{n}\to\bar{\alpha}_{ii}\left(\left\{N_j\right\}_{j=1}^{n}\right)\in\left[\alpha_{ii}^{l},\alpha_{ii}^{u}\right]\\ 1-\bar{\alpha}_{ii}\left(\left\{N_j\right\}_{j=1}^{n}\right)N_i^{l}-\sum_{\substack{j=1\\ j\neq i}}^{n}\alpha_{ij}N_j\ge 0\end{array}\right\}\\
\left\{\begin{array}{l}\forall i\in R^{-}\\ \forall\left\{N_j\right\}_{j=1}^{n}:\begin{cases}N_i=N_i^{l}\\ N_j\in\left[N_j^{l},N_j^{u}\right]\ \forall j\in\left\{1,n\right\}-\left\{i\right\}\end{cases}\end{array}\right\}\exists\bar{\alpha}_{ii}:\left\{\begin{array}{l}\left\{N_j\right\}_{j=1}^{n}\to\bar{\alpha}_{ii}\left(\left\{N_j\right\}_{j=1}^{n}\right)\in\left[\alpha_{ii}^{l},\alpha_{ii}^{u}\right]\\ 1-\bar{\alpha}_{ii}\left(\left\{N_j\right\}_{j=1}^{n}\right)N_i^{l}-\sum_{\substack{j=1\\ j\neq i}}^{n}\alpha_{ij}N_j\le 0\end{array}\right\}
\end{array}\right\}\overset{N_j^{l}>0\ \forall j=1,n}{\Leftrightarrow}
$$

$$
\left\{\begin{array}{l}
\left\{\begin{array}{l}\forall i\in R^{+}\\ \forall\left\{N_j\right\}_{j=1}^{n}:\begin{cases}N_i=N_i^{u}\\ N_j\in\left[N_j^{l},N_j^{u}\right]\ \forall j\in\left\{1,n\right\}-\left\{i\right\}\end{cases}\end{array}\right\}\left\{1-\alpha_{ii}^{u}N_i^{u}-\sum_{\substack{j=1\\ j\neq i}}^{n}\alpha_{ij}N_j\le 0\right\}\\
\left\{\begin{array}{l}\forall i\in R^{-}\\ \forall\left\{N_j\right\}_{j=1}^{n}:\begin{cases}N_i=N_i^{u}\\ N_j\in\left[N_j^{l},N_j^{u}\right]\ \forall j\in\left\{1,n\right\}-\left\{i\right\}\end{cases}\end{array}\right\}\left\{1-\alpha_{ii}^{l}N_i^{u}-\sum_{\substack{j=1\\ j\neq i}}^{n}\alpha_{ij}N_j\ge 0\right\}\\
\left\{\begin{array}{l}\forall i\in R^{+}\\ \forall\left\{N_j\right\}_{j=1}^{n}:\begin{cases}N_i=N_i^{l}\\ N_j\in\left[N_j^{l},N_j^{u}\right]\ \forall j\in\left\{1,n\right\}-\left\{i\right\}\end{cases}\end{array}\right\}\left\{1-\alpha_{ii}^{l}N_i^{l}-\sum_{\substack{j=1\\ j\neq i}}^{n}\alpha_{ij}N_j\ge 0\right\}\\
\left\{\begin{array}{l}\forall i\in R^{-}\\ \forall\left\{N_j\right\}_{j=1}^{n}:\begin{cases}N_i=N_i^{l}\\ N_j\in\left[N_j^{l},N_j^{u}\right]\ \forall j\in\left\{1,n\right\}-\left\{i\right\}\end{cases}\end{array}\right\}\left\{1-\alpha_{ii}^{u}N_i^{l}-\sum_{\substack{j=1\\ j\neq i}}^{n}\alpha_{ij}N_j\le 0\right\}
\end{array}\right\}\Leftrightarrow
$$

$$
\left\{\begin{array}{l}
\left\{0\ge\left\{\begin{array}{l}\underset{\left\{N_j\right\}_{\substack{j=1\\ j\neq i}}^{n}}{max}\left[1-\alpha_{ii}^{u}N_i^{u}-\sum_{\substack{j=1\\ j\neq i}}^{n}\alpha_{ij}N_j\right]\\ s.t.\ N_j^{l}\le N_j\le N_j^{u}\ \forall j=1,n;j\neq i\end{array}\right.\right\}\ \forall i\in R^{+},\left\{0\le\left\{\begin{array}{l}\underset{\left\{N_j\right\}_{\substack{j=1\\ j\neq i}}^{n}}{min}\left[1-\alpha_{ii}^{l}N_i^{l}-\sum_{\substack{j=1\\ j\neq i}}^{n}\alpha_{ij}N_j\right]\\ s.t.\ N_j^{l}\le N_j\le N_j^{u}\ \forall j=1,n;j\neq i\end{array}\right.\right\}\ \forall i\in R^{+}\\
\left\{0\le\left\{\begin{array}{l}\underset{\left\{N_j\right\}_{\substack{j=1\\ j\neq i}}^{n}}{min}\left[1-\alpha_{ii}^{l}N_i^{u}-\sum_{\substack{j=1\\ j\neq i}}^{n}\alpha_{ij}N_j\right]\\ s.t.\ N_j^{l}\le N_j\le N_j^{u}\ \forall j=1,n;j\neq i\end{array}\right.\right\}\ \forall i\in R^{-},\left\{0\ge\left\{\begin{array}{l}\underset{\left\{N_j\right\}_{\substack{j=1\\ j\neq i}}^{n}}{max}\left[1-\alpha_{ii}^{u}N_i^{l}-\sum_{\substack{j=1\\ j\neq i}}^{n}\alpha_{ij}N_j\right]\\ s.t.\ N_j^{l}\le N_j\le N_j^{u}\ \forall j=1,n;j\neq i\end{array}\right.\right\}\ \forall i\in R^{-}
\end{array}\right\}\Leftrightarrow
$$

$$\left\{\begin{array}{l}\left\{0 \geq 1-\alpha_{ii}^{u} N_{i}^{u}-\sum_{\substack{j \in \mathrm{A}_{i}^{+} \\ j \neq i}} \alpha_{ij} N_{j}^{l}-\sum_{\substack{j \in \mathrm{A}_{i}^{-} \\ j \neq i}} \alpha_{ij} N_{j}^{u}\right\} \ \forall i \in R^{+}, \left\{0 \leq 1-\alpha_{ii}^{l} N_{i}^{l}-\sum_{\substack{j \in \mathrm{A}_{i}^{+} \\ j \neq i}} \alpha_{ij} N_{j}^{u}-\sum_{\substack{j \in \mathrm{A}_{i}^{-} \\ j \neq i}} \alpha_{ij} N_{j}^{l}\right\} \ \forall i \in R^{+} \\ \left\{0 \leq 1-\alpha_{ii}^{l} N_{i}^{u}-\sum_{\substack{j \in \mathrm{A}_{i}^{+} \\ j \neq i}} \alpha_{ij} N_{j}^{u}-\sum_{\substack{j \in \mathrm{A}_{i}^{-} \\ j \neq i}} \alpha_{ij} N_{j}^{l}\right\} \ \forall i \in R^{-}, \left\{0 \geq 1-\alpha_{ii}^{u} N_{i}^{l}-\sum_{\substack{j \in \mathrm{A}_{i}^{+} \\ j \neq i}} \alpha_{ij} N_{j}^{l}-\sum_{\substack{j \in \mathrm{A}_{i}^{-} \\ j \neq i}} \alpha_{ij} N_{j}^{u}\right\} \ \forall i \in R^{-}\end{array}\right\} \text{O.E.}\Delta.$$